\documentclass[11pt]{article}
\usepackage{amssymb}
\usepackage{amsmath}
\usepackage{mathrsfs}
\usepackage{graphics}
\usepackage{graphicx}
\usepackage{subfigure}
\usepackage[pagewise]{lineno}
\usepackage[T1]{fontenc}
\usepackage{latexsym,amssymb,amsmath,amsfonts,amsthm}\usepackage{txfonts}
\newcommand{\be}{\begin{eqnarray}}
\newcommand{\ben}{\begin{eqnarray*}}
\newcommand{\en}{\end{eqnarray}}
\newcommand{\enn}{\end{eqnarray*}}

\newtheorem{theorem}{Theorem}[section]
\newtheorem{lemma}{Lemma}[section]
\newtheorem{prp}[theorem]{Proposition}
\newtheorem{cor}[theorem]{Corollary}
\newtheorem{dfn}{Definition}[section]

\newtheorem{remark}{Remark}

\begin{document}
\renewcommand{\theequation}{\arabic{section}.\arabic{equation}}
\begin{titlepage}
\title{\bf Exponential Mixing for 3D
Stochastic Primitive Equations
}
\author{Zhao Dong$^{\dag,\ddag}$\ \ Jianliang Zhai$^{\flat}$\ \ Rangrang Zhang$^{\sharp,}$\thanks{Corresponding author.}\\
{\small $^\dag$ RCSDS, Academy of Mathematics and Systems Science, Chinese Academy of Sciences, Beijing 100190, China.}\\
{\small $^\ddag$  School of Mathematical Sciences, University of Chinese Academy of Sciences.}\\
{\small $^\flat$ School of Mathematical Sciences, University of Science and Technology of China,}\\
 {\small Key Laboratory of Wu Wen-Tsun Mathematics, Chinese Academy of Sciences,}\\
 {\small No. 96 Jinzhai Road, Hefei, 230026, China.}\\
 {\small $^\sharp$ Department of  Mathematics,
Beijing Institute of Technology, Beijing, 100081, China.}\\
({\sf dzhao@amt.ac.cn},\ {\sf zhaijl@ustc.edu.cn},\ {\sf rrzhang@amss.ac.cn} )}
\date{}
\end{titlepage}
\maketitle


%


\noindent\textbf{Abstract}:
In this paper, we prove that weak solutions of 3D stochastic
primitive equations have exponential mixing property if the noise is sufficiently smooth and non-degenerate. With the help of uniqueness of strong solution of 3D stochastic
primitive equations, we obtain that all weak solutions which are limitations of Galerkin approximations share the same invariant measure. In particular, the invariant measure of strong solution is unique. The coupling method plays a key role.





\section{Introduction }
This paper concerns the ergodic theory of 3D stochastic primitive equations. The primitive equations can well model the large-scale motion of the ocean. They are derived from the Navier-Stokes equations with rotation, coupled with thermodynamics
and salinity diffusion-transport equations.
Beyond their considerable significance in physical applications,
the primitive equations have generated much interest from the mathematics community due to their rich nonlinear, nonlocal character and their anisotropic structure. For example, the mathematical study of the primitive equations originated in a series of articles by Lions,
Temam and Wang in the early 1990s (see \cite{L-T-W-1,L-T-W-2,L-T-W-3,L-T-W-4} and references therein), where they set up the mathematical framework and showed the global existence of weak solutions.
Cao and Titi \cite{C-T-1} developed an approach to deal with the $L^6$-norm of the fluctuation $\tilde{v}$ of horizontal velocity and obtained the global well-posedness of the 3D viscous primitive equations.


The addition of a stochastic noise to this physical model is fully natural as it represents external random perturbations or a lack of knowledge of certain physical parameters.
Along with the great successful developments of deterministic
primitive equations, the random situation has also been developed rapidly.
Guo and Huang \cite{Guo} obtained
the existence of universal random attractor of strong solution under the assumptions that that the momentum equation is driven by an additive stochastic forcing and the
thermodynamical equation is driven by a fixed heat source.
Debussche, Glatt-Holtz, Temam and Ziane \cite{D-G-T-Z} established the global well-posedness of strong solution to 3D stochastic primitive equations. Dong, Zhai and Zhang \cite{DZZ} established the Freidlin-Wentzell's large deviations for 3D stochastic primitive equations.

The purpose of this paper is to prove the exponential mixing property of weak solutions of 3D stochastic primitive equations.
There are several works about the ergodicity of 3D stochastic primitive equations. We mention two of them which are relevant to our paper.
Tachim Medjo \cite{M-1} proved that weak solutions of 3D stochastic primitive equations converges to the stationary solution exponentially in the mean square under the conditions that the viscosity is large enough (bigger than $L^2-$norm of the stationary solution) and the covariance operator of the noise satisfies some exponential decay property. Glatt-Holtz, Kukavica, Vicol and Ziane \cite{GKVZ} established the existence of invariant measures of the strong solution of 3D stochastic primitive equations. In this paper, we obtain two new results as follows:
\begin{description}
  \item[1.] The exponential mixing property holds for special weak solutions to 3D stochastic primitive equations. The viscosity is only assumed to be strictly positive, which relaxes the condition required by  \cite{M-1}. Moreover, all weak solutions which are limits of Galerkin approximations share the same stationary measure.
  \item[2.] The invariant probability measure of the strong solution to 3D stochastic primitive equations is unique, which supplements \cite{GKVZ}.
\end{description}

We adopt the coupling method introduced by Odasso \cite{O-C} to establish the exponential ergodicity. The key idea of the coupling method is to construct a coupling of Galerkin approximating solutions with small initial data in an appropriate norm, which has the property that the probability of their meeting has a uniform lower bound $\frac{3}{4}$ (see Proposition \ref{prp1}). Compared with \cite{O-C},
we need to choose a small ball in $\mathbb{H}_3$ instead of $\mathbb{H}_2$ because the nonlinear term $-\int^{z}_{-1}\nabla_H\cdot v(t,x,y,z')dz'\frac{\partial v}{\partial z}$ of 3D stochastic primitive equations is one order higher than $v\frac{\partial v}{\partial z}$ of 3D stochastic Navier-Stokes equations. Furthermore, we have to show that the exist time into this ball admits an exponential moment. During the proof process, we have to deal with some high order Sobolev norm estimates, such as  $\|\cdot\|^6,\ \|\cdot\|_2^4 $ and $\ \|\cdot\|_3^4$. These are highly non-trivial.
To prove the uniqueness of invariant measures, the uniqueness of the strong solution to 3D stochastic primitive equations plays an important role. It's worth mention that such result has not been established for 3D  stochastic Navier-Stokes equations.

\vskip0.3cm
This paper is organized as follows.
In Sect. 2, we introduce the mathematical framework of 3D stochastic primitive equations and main results.
In Sect. 3, we give the proof of weak solution to 3D stochastic primitive equations.
In Sect. 4, the exponential mixing property for Galerkin sequences from $\mathbb{H}_3$ are proved.
In Sect. 5, we show the uniqueness of stationary probability measure.
The proof process of two propositions in Sect. 4 is presented in Sect. 6.


\section{Preliminaries and the statement of main results}
The 3D stochastic primitive equations under a stochastic forcing, in a Cartesian system, are written as
\begin{eqnarray}\label{eq-1}
\frac{\partial v}{\partial t}+(v\cdot \nabla_H)v+\theta\frac{\partial v}{\partial z}+f{k}\times v +\nabla_H P - \nu_1 \Delta v &=&\sigma_1(v,T)\frac{dW_1}{dt},\\
\label{eq-2}
\partial_{z}P+T&=&0,\\
\label{eq-3}
\nabla_H\cdot v+\partial_{z}\theta&=&0,\\
\label{eq-4}
\frac{\partial T}{\partial t}+(v\cdot\nabla_H)T+\theta\frac{\partial T}{\partial z}- \nu_2\Delta T&=&\sigma_2(v,T)\frac{dW_2}{dt},
\end{eqnarray}
where the horizontal velocity field $v=(v_{1},v_{2})$, the velocity $\theta$, the temperature\ $T$ and the pressure\ $P$ are unknown functions. $f$ is the Coriolis parameter. ${k}$ is vertical unit vector. $\nu_1$ and $\nu_2$ are the viscosity. $\nabla_H=(\partial_x,\partial_y)$  to be the horizontal gradient operator and $\Delta=\partial^{2}_{x}+\partial^{2}_{y}+\partial^{2}_{z}$ to be the three dimensional Laplacian. $W_1$ and $W_2$ are two independent cylindrical Wiener processes on $H_1$ and $H_2$, respectively. The spaces $H_1$ and $H_2$ will be defined below. The spatial variable $(x,y,z)$ belongs to $\mathcal{M}:= \mathbb{T}^2\times (-1,0)$.

%
As in \cite{C-L-T}, the boundary and initial conditions of $(\ref{eq-1})-(\ref{eq-4})$ are given by
\begin{eqnarray}\label{eqq1}
v, \ \theta\ {\rm{and}} \ T\ {\rm{are\ periodic\ in}}\ x\ {\rm{and}}\ y,\\
\label{eqq2}
(\partial_{z}v,\theta)\mid _{z=-1,0}=(0,0),\ T\mid_{z=-1}=1, \ T\mid_{z=0}=0.\\
\label{eqq3}
(v,T)\mid_{t=0}=(v_0,T_0).
\end{eqnarray}
Replacing $T$ and $P$ by $T+z$ and $P-\frac{z^2}{2}$, respectively, then $(\ref{eq-1})-(\ref{eqq3})$ is equivalent to the following system
\begin{eqnarray}\label{eqq4}
\frac{\partial v}{\partial t}+(v\cdot \nabla_H)v+\theta\frac{\partial v}{\partial z}+f{k}\times v +\nabla_H P -\nu_1\Delta v &=&\sigma_1(v,T+z)\frac{dW_1}{dt},\\
\label{eqq5}
\partial_{z}P+T&=&0,\\
\label{eqq6}
\nabla_H\cdot v+\partial_{z}\theta&=&0,\\
\label{eqq7}
\frac{\partial T}{\partial t}+(v\cdot\nabla_H)T+\theta(\frac{\partial T}{\partial z}+1)-\nu_2\Delta T&=&\sigma_2(v,T+z)\frac{dW_2}{dt},
\end{eqnarray}
subject to the boundary and initial conditions
\begin{eqnarray}\label{eqq8}
v, \ \theta\ {\rm{and}} \ T\ {\rm{are\ periodic\ in}}\ x\ {\rm{and}}\ y,\\
\label{eqq9}
(\partial_{z}v,\theta)\mid _{z=-1,0}=(0,0),\ T\mid_{z=-1,z=0}=0,\\
 \label{eqq10}
(v,T)\mid_{t=0}=(v_0,T_0).
\end{eqnarray}
We still denote by $T_0$ the initial temperature in (\ref{eqq10}), although it is now different from that in $(\ref{eqq3})$.

Inherent symmetries in the equations show that the solution of the primitive equations on $\mathbb{T}^2\times (-1,0)$ with boundaries (\ref{eqq8})-(\ref{eqq10}) can be recovered by solving the equations with periodic boundary conditions in $x, y$ and $z$ variables on the extended domain $\mathbb{T}^2\times (-1,1):= \mathbb{T}^3$ and restricting to $z\in (-1,0)$.

To see this, consider any solution of (\ref{eqq4})-(\ref{eqq7}) with boundaries (\ref{eqq8})-(\ref{eqq10}), we perform that
\begin{eqnarray*}
v(x,y,z)&=&v(x,y,-z), \ \rm{for}\ (x,y,z)\in \mathbb{T}^2\times (0,1),\\
T(x,y,z)&=&-T(x,y,-z), \ \rm{for}\ (x,y,z)\in \mathbb{T}^2\times (0,1),\\
P(x,y,z)&=&P(x,y,-z), \ \rm{for}\ (x,y,z)\in \mathbb{T}^2\times (0,1),\\
\theta(x,y,z)&=&-\theta(x,y,-z),\ \rm{for}\ (x,y,z)\in \mathbb{T}^2\times (0,1).
\end{eqnarray*}
 We also extend $\sigma_1$ in the even fashion and  $\sigma_2$ in the odd fashion across $\mathbb{T}^2\times \{0\}$. Hence, we consider the primitive equations on the extended domain $\mathbb{T}^3=\mathbb{T}^2\times (-1,1)$,
 \begin{eqnarray}\label{eqq14}
\frac{\partial v}{\partial t}+(v\cdot \nabla_H)v+\theta\frac{\partial v}{\partial z}+f{k}\times v +\nabla_H P -\nu_1\Delta v &=&\sigma_1({v},{T})\frac{dW_1}{dt},\\
\label{eqq15}
\partial_{z}P+T&=&0,\\
\label{eqq16}
\nabla_H\cdot v+\partial_{z}\theta&=&0,\\
\label{eqq17}
\frac{\partial T}{\partial t}+(v\cdot\nabla_H)T+\theta(\frac{\partial T}{\partial z}+1)-\nu_2\Delta T&=&\sigma_2({v},{T}+z)\frac{dW_2}{dt},
\end{eqnarray}
 subject to the boundary and initial conditions
 \begin{eqnarray}\label{eqq11}
v, \ \theta\ , P\ {\rm{and}}\ T\ {\rm{are\ periodic\ in}}\  x\ ,\ y,\ z,\\
 \label{eqq12}
v\ {\rm{and}}\ P\ {\rm{are\ even\ in}} \ z, \  \theta\ {\rm{and}}\ T\ {\rm{are\ odd\ in}}\ z,\\
\label{eqq13}
(v,T)\mid_{t=0}=(v_0,T_0),
\end{eqnarray}
where $T_0$ denotes the initial temperature in (\ref{eqq10}), although it is now different from that in $(\ref{eqq3})$.

Due to the equivalence of the above two kinds boundary conditions, we consider the system (\ref{eqq14})-(\ref{eqq13}) throughout the whole paper. Notice that (\ref{eqq12}) is a symmetry condition, which is preserved by system (\ref{eqq14})-(\ref{eqq17}), i.e., if a smooth solution to the system (\ref{eqq14})-(\ref{eqq17}) exists and is unique, then it satisfies the symmetry condition (\ref{eqq12}) as long as it is initially satisfied.
\subsection{Reformulation}\label{se-2}
With the help of the incompressibility condition (\ref{eqq16}) and the symmetry condition (\ref{eqq12}), the vertical velocity $\theta$ can be expressed in terms of the horizonal velocity $v$ as
\begin{equation}
\theta(t,x,y,z)=\Phi(v)(t,x,y,z)=-\int^{z}_{-1}\nabla_H\cdot v(t,x,y,z')dz'.
\end{equation}
Moreover,
\[
\int^{1}_{-1}\nabla_H\cdot v  dz=0.
\]
Supposing that $p_{b}$ is a certain unknown function at $\Gamma_{b}$, and integrating (\ref{eqq15}) from $-1$ to $z$, we have
\[
P(x,y,z,t)= p_{b}(x,y,t)-\int^{z}_{-1} T(x,y,z',t) dz'.
\]
Then, we make scaling transformation of $T$. Let $S=\sqrt{C_0} T$, where $C_0$ is a positive constant will be given in Section \ref{l-2}. In this case, $(\ref{eqq14})-(\ref{eqq13})$ can be rewritten as
 \begin{eqnarray}\label{eq5-1}
&\frac{\partial v}{\partial t}+(v\cdot \nabla_H)v+\Phi(v)\frac{\partial v}{\partial z}+f{k}\times v +\nabla_H p_{b}-\frac{1}{\sqrt{C_0}}\int^{z}_{-1}\nabla_H S dz' -\nu_1\Delta v =\phi(v,S),&\\
\label{eq-6-1}
&\frac{\partial S}{\partial t}+(v\cdot\nabla_H)S+\Phi(v)\frac{\partial S}{\partial z}+ \Phi(v) -\nu_2\Delta S=\varphi(v,S),&\\
\label{eq-7-1}
&\int^{1}_{-1}\nabla_H\cdot v  dz=0,&
\end{eqnarray}
where
\[
\phi(v,S)=\sigma_1({v},\frac{1}{\sqrt{C_0}}{S})\frac{dW_1}{dt}, \quad \varphi(v,S)=\sqrt{C_0}\sigma_2({v},\frac{1}{\sqrt{C_0}}{S}+z)\frac{dW_2}{dt}.
\]
The boundary and initial conditions for $(\ref{eq5-1})-(\ref{eq-7-1})$ are given by
\begin{eqnarray}\label{eqq18}
v\ {\rm{and}}\ S\ {\rm{are\ periodic\ in}}\ x, \ y \ {\rm{and}} \ z,\\
\label{eqq19}
v\ {\rm{and}}\ P\ {\rm{are\ even\ in}} \ z, \  \theta\ {\rm{and}}\ T\ {\rm{are\ odd\ in}}\ z,\\
\label{eqq20}
(v, S)\mid_{t=0}=(v_0,S_0).
\end{eqnarray}
To prove the exponential mixing property of  $(v,T)$ of the system $(\ref{eqq14})-(\ref{eqq13})$, it suffices to prove $(v, S)$ of $(\ref{eq5-1})-(\ref{eqq20})$ has the exponential mixing property.
\subsection{Working spaces}
 Let $\mathcal{L}(K_1;K_2)$ (resp. $\mathcal{L}_2(K_1;K_2)$) be the space of bounded (resp. Hilbert-Schmidt) linear operators from the Hilbert space $K_1$ to $K_2$. $|\cdot|_{L^2(\mathbb{T}^2)}$ and $|\cdot|_{H^p(\mathbb{T}^2)}$ stand for the norms of $L^2(\mathbb{T}^2)$ and $H^p(\mathbb{T}^2)$, respectively. Let $|\cdot|_{p}$ be the norm of $L^p(\mathbb{T}^3)$ for integer number $p\in (0,\infty)$. In particular, $|\cdot|$ and $(\cdot,\cdot)$ represent the norm and inner product of $L^2(\mathbb{T}^3)$. For the classical Sobolev space $H^p(\mathbb{T}^3)=W^{p,2}(\mathbb{T}^3)$, $p\in \mathbb{N}_+$,
\begin{equation}\notag
\left\{
  \begin{array}{ll}
    H^p(\mathbb{T}^3)=\Big\{Y\in L^2(\mathbb{T}^3)\Big| \partial_{\alpha}Y\in L^2(\mathbb{T}^3)\ {\rm for} \ |\alpha|\leq p\Big\},&  \\
    |Y|^2_{H^p(\mathbb{T}^3)}=\sum_{0\leq|\alpha|\leq p}|\partial_{\alpha}Y|^2. &
  \end{array}
\right.
\end{equation}
It's  known that $(H^p(\mathbb{T}^3), |\cdot|_{H^p(\mathbb{T}^3)})$ is a Hilbert space.

Define working spaces for the system $(\ref{eq5-1})-(\ref{eqq20})$ as
 \begin{eqnarray}\notag
 &&\mathcal{V}_1:=\left\{v\in (C^{\infty}(\mathbb{T}^3))^2;\ \int^1_{-1}\nabla_H\cdot v dz=0, v \ {\rm{is\ periodic\ in}}\ x, \ y\ {\rm{and \ even\ in}} \ z, \int_{\mathbb{T}^3}vdxdydz=0\right\},\\ \notag
 &&\mathcal{V}_2:=\left\{S\in C^{\infty}(\mathbb{T}^3);\  \ S\ {\rm{is\ periodic\ in}} \ x, \ y \ {\rm{and\ odd\ in}} \ z, \int_{\mathbb{T}^3}Sdxdydz=0 \right\},
 \end{eqnarray}
$V_1$= the closure of $ \mathcal{V}_1$ with respect to the norm $|\cdot|_{H^1(\mathbb{T}^3)}\times |\cdot|_{H^1(\mathbb{T}^3)}$,\\
$V_2$= the closure of $ \mathcal{V}_2$ with respect to the norm $|\cdot|_{H^1(\mathbb{T}^3)}$,\\
$H_1$= the closure of $ \mathcal{V}_1$ with respect to the norm $|\cdot|\times |\cdot|$,\\
$H_2$= the closure of $ \mathcal{V}_2$ with respect to the norm $|\cdot|$,\\
\[
V=V_1\times V_2, \quad H=H_1\times H_2.
\]

 The inner products and norms on $V$, $H$ are given by
\begin{eqnarray*}
(Y,\tilde{Y})_{V}&=&(v,\tilde{v})_{V_1}+(S,\tilde{S})_{V_2},\\
(Y,\tilde{Y})&=&(v,\tilde{v})+(S,\tilde{S})=(v_1, \tilde{v}_1)+(v_2, \tilde{v}_2)+(S,\tilde{S}),\\
(Y,Y)^{\frac{1}{2}}_{V}&=&(v,v)^{\frac{1}{2}}_{V_1}+(S,S)^{\frac{1}{2}}_{V_2}, \quad \|Y\|_{V}=(Y,Y)^{\frac{1}{2}}_V.
\end{eqnarray*}
where $Y=(v,S), \tilde{Y}=(\tilde{v},\tilde{S})$. By the Riesz representation theorem, we can identify the dual space $H'$ of $H$, which implies that
\[
V\subset H = H'\subset V',
\]
where the two inclusions are compact continuous.
\subsection{Some Functionals}
Define three bilinear forms $a:V\times V\rightarrow \mathbb{R}$,\ $a_1:V_1\times V_1\rightarrow \mathbb{R}$,\ $a_2:V_2\times V_2\rightarrow \mathbb{R}$,
 and their corresponding linear operators $A: V\rightarrow V^{'}$, $A_1: V_1\rightarrow V^{'}_1$, $A_2: V_2\rightarrow V^{'}_2$ by setting
 \[
 a(Y,Y_1):=(AY,Y_1)=   
   a_1(v,v_1)+  
  a_2(S,S_1)\\  
  ,
 \]
 where
\begin{eqnarray}\notag
a_1(v,v_1):=(A_1v, v_1)=\nu_1\int_{\mathbb{T}^3}\nabla v\cdot \nabla v_1dxdydz,
\end{eqnarray}
\begin{eqnarray}\notag
a_2(S,S_1):=(A_2S, S_1)=\nu_2\int_{\mathbb{T}^3}\nabla S\cdot \nabla S_1dxdydz,
\end{eqnarray}
for any $Y=(v,S)$, $Y_1=(v_1, S_1)\in V$.

It's known that $A_1$ is a self-adjoint operator with discrete spectrum in $H_1$. Denote by $\{k_n\}_{n=1,2,\cdot\cdot\cdot}$ the eigenbasis of $A_1$ and suppose its associated eigenvalues $\{\gamma_n\}_{n=1,2,\cdot\cdot\cdot}$ is increasing. Similarly, $A_2$ is a self-adjoint operator with discrete spectrum in $H_2$. Let $(l_n)_{n=1,2,\cdot\cdot\cdot}$ be the eigenbasis of $A_2$ with increasing corresponding eigenvalues $\{\lambda_n\}_{n=1,2,\cdot\cdot\cdot}$. Denote
 $\bar{e}_{n,0}=\left(                 
  \begin{array}{c}   
   k_n\\  
   0\\  
  \end{array}
\right)$ and $\bar{e}_{0,m}=\left(                 
  \begin{array}{c}   
   0\\  
   l_m\\  
  \end{array}
\right)$,
it is easy to know $\{\bar{e}_{n,0},\bar{e}_{0,m}\}_{n,m=1,2,\cdots}$ is an eigenbasis of $(A, {D}(A))$.
By means of rearrangement, we can construct an eigenbasis of $(A, {D}(A))$ denoted $\{e_n\}_{n=1,2,\cdots}$ such that the associated eigenvalues $\{\mu_n\}_{n=1,2,\cdot\cdot\cdot}$ is an increasing sequence.

For any $s\in \mathbb{R}$, the fractional power $(A^s, {D}(A^s))$ of the operator  $(A,{D}(A) )$ is defined as
\begin{equation}\notag
\left\{
  \begin{array}{ll}
    {D}(A^s)=\Big\{Y=\sum_{n=1}^{\infty}y_n e_n \Big|  \sum_{n=1}^{\infty}\mu^{2s}_n|y_n|^2<\infty\Big\}; &  \\
    A^s Y=\sum_{n=1}^{\infty}\mu^s_ny_n e_n, \quad where\  Y=\sum_{n=1}^{\infty}y_n e_n .&
  \end{array}
\right.
\end{equation}
Set
\[
\|Y\|^{A}_s=|A^{\frac{s}{2}}Y|,\quad \mathbb{H}^A_s={D}(A^{\frac{s}{2}}),
\]
then $(\mathbb{H}^A_0, \|\cdot\|^A_0)=(H,|\cdot|)$ and $(\mathbb{H}^A_1, \|\cdot\|^A_1)=(V, \|\cdot\|_V)$. For simplicity, denote $\|\cdot\|=\|\cdot\|_V$. It's obvious that $(\mathbb{H}^A_s, \|\cdot\|^A_s)$ is a Hilbert space.
Similarly, we can define $(\mathbb{H}^{A_1}_s,\|\cdot\|^{A_1}_s)$ and $(\mathbb{H}^{A_2}_s, \|\cdot\|^{A_2}_s)$. For convenience,
all of them will be denoted by $(\mathbb{H}_s,\|\cdot\|_s)$ for $s\in \mathbb{R}$.

Now, we define three functionals $b: V\times V\times V\rightarrow \mathbb{R}$, $b_i: V_1\times V_i\times V_i\rightarrow \mathbb{R}\ (i=1,2)$  and the associated operators $B: V\times V\rightarrow V'$, $B_i: V_1\times V_i\rightarrow V'_i\ (i=1,2)$ by setting
\begin{eqnarray*}
b(Y,Y_1,Y_2)&:=&(B(Y,Y_1),Y_2)=   
   b_1(v,v_1,v_2)+  
   b_2(v, S_1, S_2) 
 ,\\
b_1(v,v_1,v_2)&:=&(B_1(v,v_1), v_2)=\int_{\mathbb{T}^3}\left[(v\cdot \nabla_H)v_1+\Phi(v)\frac{\partial v_1}{\partial z}\right]\cdot v_2dxdydz,\\
b_2(v, S_1, S_2)&:=&(B_2(v,S_1), S_2)=\int_{\mathbb{T}^3}\left[(v\cdot \nabla_H)S_1+\Phi(v)\frac{\partial S_1}{\partial z}\right] S_2dxdydz,
\end{eqnarray*}
for any $Y=(v, S)$, $Y_i=(v_i,S_i)\in V$. Then we have

\begin{lemma}\label{lemma-1}
 For any $Y\in V$, $Y_1\in V$,
\[
\left(B(Y,Y_1),Y_1\right)=b(Y,Y_1,Y_1)=b_1(v,v_1,v_1)=b_2(v,S_1,S_1)=0.
\]
\end{lemma}
Moreover, we define another functional $g: V\times V \rightarrow \mathbb{R}$ and the associated linear operator $G: V\rightarrow V'$ by
\begin{eqnarray*}\notag
g(Y,Y_1)&:=&(G(Y), Y_1)\notag
\\      &=&\int_{\mathbb{T}^3}\left[f(k\times v)\cdot v_1+\left(\nabla_H p_b-\frac{1}{\sqrt{C_0}}\int^z_{-1}\nabla_H Sdz'\right)\cdot v_1+\Phi(v)\cdot S_1\right]dxdydz 
.
\end{eqnarray*}
We deduce from (\ref{eq-7-1}) that
\begin{equation}\notag
(v,\nabla_H p_{b})=\left(\int^{1}_{-1}vdz,\nabla_H p_b\right)_{L^2(\mathbb{T}^2)}=-\left(p_b,\int^{1}_{-1}\nabla_H \cdot vdz \right)_{L^2(\mathbb{T}^2)}=0.
\end{equation}
Since $(v,f{k}\times v)=0$, we have
\begin{lemma}\label{lemma-2}
\begin{description}\notag
  \item[(i)]
\[
g(Y,Y)=(G(Y),Y)=
-\frac{1}{\sqrt{C_0}}\int_{\mathbb{T}^3}\Big[\left(\int^z_{-1}\nabla_H Sdz'\right)\cdot v+\Phi(v)\cdot S\Big]dxdydz
.
\]
  \item[(ii)] For $Y=(v,S), \tilde{Y}=(\tilde{v}, \tilde{S})$, there exists a constant $C$ such that
\begin{eqnarray}\label{a-2}
|(G(Y),Y)|&\leq& C(|S|\|v\|\vee\|S\||v|),\\
\label{a-3}
|(G(Y),\tilde{Y})|&\leq& C|v||\tilde{v}|+C(|S|\|\tilde{v}\|\vee\|S\||\tilde{S}|).
\end{eqnarray}
\end{description}
\end{lemma}
Using the above functionals, (\ref{eq5-1}) and (\ref{eq-6-1}) can be written as
\begin{eqnarray}\label{aa}
\left\{
  \begin{array}{ll}
    dY(t)+AY(t)dt+B(Y(t),Y(t))dt+G(Y(t))dt=\Psi(Y(t))dW(t), \\
    Y(0)=y_0,
  \end{array}
\right.
\end{eqnarray}
where
\begin{equation}\notag
W=\left(                 
  \begin{array}{c}   
    W_1\\  
    W_2 \\  
  \end{array}
\right) ,\quad
\Psi(Y)
=\left(                 
  \begin{array}{cc}   
   \phi(v,S) & 0 \\  
   0 & \varphi(v,S) \\  
  \end{array}
\right).
\end{equation}





\subsection{Hypotheses}\label{l-2}
Recall that $W$ is said to be a $(\mathcal{F}_t)_{t\in [ 0,T]}-$cylindrical Wiener process on $H$ if $W$ is $(\mathcal{F}_t)_t-$adapted, $W(\cdot+t)-W(t)$ is independent of $\mathcal{F}_t$ for any $t\geq 0$ and $W$ is a cylindrical Wiener process on $H$. Let $E$ be a Polish space.

Suppose $C_0\geq \frac{8}{\lambda_1}$ in Sect. \ref{se-2}.  $W_1$, $W_2$ are two independent cylindrical Wiener processes on $H_1$ and $H_2$, respectively. $W=(W_1, W_2)^\perp$ can be written as  $W=\sum^{\infty}_{n=1}\beta_n e_n$, where $\{\beta_n\}_n$ is a sequence of 1d real-valued standard Brownian motions. The covariance operator $\Psi$  satisfies
\begin{description}
   \item[\textbf{Hypothesis H0}] $\Psi:$ $H\rightarrow \mathcal{L}_2(H;H)$ is a continuous and bounded Lipschitz mapping, i.e.,
       \[
       \|\Psi(y)\|^2_{\mathcal{L}_2(H;H)}\leq \lambda_0|y|^2 +\rho,\ y\in H,
       \]
for some constants $\lambda_0\geq0,\ \rho\geq0.$
\end{description}
 \begin{description}
   \item[\textbf{Hypothesis H1}]
        \begin{description}
          \item[(i)]  There exists $\varepsilon_0>0$ and a family $\{\Psi_n\}_{n=1,2,\cdot\cdot\cdot}$ of continuous mappings $H\rightarrow \mathbb{R}$ with continuous Fr\'{e}chet  derivatives such that
\[
\left\{
  \begin{aligned}
    \Psi(y)dW&=\sum^{\infty}_{n=1}\Psi_n(y)e_nd\beta_n \quad {\rm{where}} \quad  W=\sum^{\infty}_{n=1}\beta_n e_n, \\
    \kappa_0&=\sum^{\infty}_{n=1}\sup_{y\in H}|\Psi_n(y)|^2\mu^{2+\varepsilon_0}_n<\infty.
  \end{aligned}
\right.
\]
\item[(ii)]  There exists $\kappa_1$ such that for any $y, \eta \in \mathbb{H}_3$,
\[
\sum^{\infty}_{n=1}|\Psi'_n(y)\cdot \eta|^2\mu^3_n<\kappa_1\|\eta\|^2_3.
\]
          \item[(iii)] For any $y\in H$ and $n\in \mathbb{N}$,
\[
\Psi_n(y)>0, \quad \kappa_2=\sup_{y\in H}\|\Psi^{-1}(y)\|^2_{\mathcal{L}(\mathbb{H}_4; H)}< \infty,
\]
where
\[
\Psi^{-1}(y)\cdot h=\sum^{\infty}_{n=1}\Psi_n^{-1}(y)h_n e_n \quad {\rm for} \quad h=\sum^{\infty}_{n=1}h_n e_n.
\]
        \end{description}

Set $\kappa=\kappa_0 +\kappa_1+\kappa_2+1$.
\end{description}
\begin{remark}  \textbf{Hypothesis H1} implies \textbf{Hypothesis H0} and $\Psi=A^{-\frac{\beta}{2}}$ fulfills \textbf{Hypothesis H1}, provided $\beta\in (\frac{7}{2},4]$.
\end{remark}
\subsection{Main Results}
The main results in this paper are as follows.
\begin{theorem}\label{thm-1}[\textbf{Existence of weak solutions}]
Let the initial data $y_0\in H$. Assume \textbf{Hypothesis H0} is in force, there exists a weak solution of (\ref{aa}) in the sense of Definition \ref{dfn1}.
\end{theorem}
Let $\lambda$ and $Y$ be probability measure and a random variable on $(\Omega, \mathcal{F})$, respectively. We denote by $P(H)$ the set of probability measure on $H$ endowed with the Borelian $\sigma-$algebra. $\mathcal{D}_{\mathbb{P}_{\lambda}}(Y)$ denotes the law of $Y$ under $\mathbb{P}_{\lambda}$.  If $\mu= \mathcal{D}_{\mathbb{P}_{\mu}}(Y(t))$, for any $t\geq 0$, then $\mu\in P(H)$ is said to be a stationary measure.
\begin{theorem}\label{th-1}[\textbf{Exponential mixing}]
 Assume \textbf{Hypothesis H1} holds. There exist positive constants $C=C(\nu_1, \nu_2, \Psi, \mathbb{T}^3)$ and $\gamma=\gamma(\nu_1, \nu_2, \Psi, \mathbb{T}^3)$ such that, for any initial law $\lambda \in {P}(H)$ and weak solution $Y(t)$ which is a limit of Galerkin approximations of (\ref{aa-1}), there exists a unique weak stationary solution $\mathbb{P}_{\mu}$ with the initial law $\mu\in {P}(H)$, such that
 \begin{eqnarray}\label{eq-12}
 \|\mathcal{D}_{\mathbb{P}_{\lambda}}(Y(t))-\mu\|_{var}\leq C{\rm{e}}^{-\gamma t}\left(1+\int_{H}|y|^2\lambda(dy)\right),
 \end{eqnarray}
provided
\[
\int_{H}|y|^2\lambda(dy)<\infty,
\]
where  $\|\cdot\|_{var}$ is the total variation norm associated to the space $\mathbb{H}_s$ for $s<-3$.
\end{theorem}

\begin{cor}\label{cor-1}
All weak solutions which are limits of Galerkin approximations share the same stationary measure. That is, $\mu$ in Theorem \ref{th-1} is independent of $\lambda$ and $\{N_k\}$. In particular, the invariant probability measure of the strong solution of (\ref{aa}) is unique.
\end{cor}

\begin{remark} Our result is not influenced by the size of the viscosity $\nu_1$ and $\nu_2$, which are only required to be strictly  positive. For simplicity, we assume $\nu_1=\nu_2=1$. Moreover,
the result of Corollary \ref{cor-1} depends on the uniqueness of the strong solution to 3D stochastic primitive equations.
\end{remark}
\subsection{Some Inequalities}
The following lemma (see \cite{Resnick}, Lemma A.4) plays an important role.
\begin{lemma}\label{lemm-2}
Let $\Lambda=(-\Delta)^{\frac{1}{2}}$. Suppose that $s>0$ and $p\in(1,\infty)$. If $f,g\in C^{\infty}(\mathbb{T}^3)$, then
\[
|\Lambda^s(fg)|_p\leq C(|f|_{p_1}|\Lambda^sg|_{p_2}+|g|_{p_3}|\Lambda^sf|_{p_4}),
\]
with $p_i\in(1,\infty]$, $i=1,\cdot\cdot\cdot,4$ such that
\[
\frac{1}{p}=\frac{1}{p_1}+\frac{1}{p_2}=\frac{1}{p_3}+\frac{1}{p_4}.
\]
\end{lemma}
Refer to \cite{C-T-1}, we have
\begin{lemma}\label{le-2}
 If $v_1 \in H^1(\mathbb{T}^3),v_2 \in H^3(\mathbb{T}^3),v_3 \in H^3(\mathbb{T}^3)$, then
\begin{description}
 \item[(i)] $|\int_{\mathbb{T}^3}v_3\cdot[(v_1 \cdot \nabla_H)v_2]dxdydz|\leq c|\nabla v_2||v_3|_3|v_1|_6\leq c |\nabla v_2||v_3|^{\frac{1}{2}}|\nabla v_3|^{\frac{1}{2}}|\nabla v_1|$,
 \item[(ii)] $|\int_{\mathbb{T}^3}\Phi(v_1)v_{2z}\cdot v_3 dxdydz|\leq c|\nabla v_1||v_3|^{\frac{1}{2}}|\nabla v_3|^{\frac{1}{2}}|\partial_{z} v_2|^{\frac{1}{2}}|\nabla \partial_{z} v_2|^{\frac{1}{2}}$.
\end{description}
\end{lemma}

\section{ Existence of weak solutions}\label{section-1}
\begin{dfn}(weak solutions).\label{dfn1}
We say that there exists a weak solution of (\ref{aa}) if for any initial law $\lambda\in P(H)$ and $T>0$, there exists a stochastic basis
 $(\Omega, \mathcal{F},\{\mathcal{F}_t\}_{t\in [ 0,T]}, \mathbb{P}_{\lambda})$, a $(\mathcal{F}_t)_t-$cylindrical Wiener process $W$ on $H$ under $\mathbb{P}_{\lambda}$ and a progressively measurable process $Y: [0,T]\times \Omega\rightarrow H$  such that
 \begin{description}
   \item[(i)] the law of $Y(0)$ under $\mathbb{P}_{\lambda}$ is $\lambda$,
   \item[(ii)]
$ Y\in L^{\infty}([0,T]; H)\cap L^2([0,T]; V)\cap C\left([0,T]; D(A^{-\frac{\alpha}{2}})\right),$
 where $\alpha > {3}$ is any fixed positive number.
   \item[(iii)] for any $t\in[0,T]$ and $\psi \in D(A^{\frac{\alpha}{2}})$, the following holds $\mathbb{P}_{\lambda}-$a.s.
\begin{eqnarray*}
(Y(t), \psi)&+&\int^{t}_{0}\Big(Y(s),A\psi\Big)ds+\int^{t}_{0}\Big(B(Y(s),Y(s)),\psi\Big)ds+\int^{t}_{0}\Big(G(Y(s)),\psi\Big)ds\notag
\\ &  & =  (Y(0), \psi)+\int^{t}_{0}\Big(\Psi(Y(s))dW(s),\psi\Big).
\end{eqnarray*}
 \end{description}
 When the initial law $\lambda$ is not specified, $y_0$ is the initial value of the weak solution $\mathbb{P}_{y_0}$.
\end{dfn}
\begin{remark}
These solutions are weak in both probability and PDE sense.

\end{remark}
To prove the existence of weak solutions of (\ref{aa}), some compact embedding results of certain functional spaces are needed (see \cite{F-G}).

Let $H$ be a separable Hilbert space with norm $|\cdot|_H$.
Given $p>1, \alpha\in (0,1)$,  let $W^{\alpha,p}([0,T]; H)$ be the Sobolev space of all $u\in L^p([0,T];H)$ such that
\[
\int^T_0\int^T_0\frac{|u(t)-u(s)|^{ p}_H}{|t-s|^{1+\alpha p}}dtds< \infty,
\]
endowed with the norm
\[
\|u\|^p_{W^{\alpha,p}([0,T]; H)}=\int^T_0|u(t)|^p_Hdt+\int^T_0\int^T_0\frac{|u(t)-u(s)|^{ p}_H}{|t-s|^{1+\alpha p}}dtds.
\]

Let $(\Omega, \mathcal{F}, \{\mathcal{F}_t\}_{t\geq 0}, \mathbb{P})$ be a stochastic basis (with expectation $\mathbb{E}$). $K$ is a separable Hilbert space and $W$ is a cylindrical Wiener process with values in $K$ defined on this stochastic basis.
For any progressively measurable process $f\in L^2(\Omega\times [0,T]; \mathcal{L}_2(K;H))$, we denote by $I(f)$ the It\^{o} integral
\[
I(f)(t)=\int^{t}_{0}f(s)dW(s) \quad t\in [0,T].
\]
Clearly, $I(f)$ is a progressively measurable process in $L^2(\Omega\times [0,T]; H)$.
\begin{lemma}\label{lemma-5}
Given $p\geq 2, \alpha<\frac{1}{2}$. Then for any progressively measurable process $f\in L^p(\Omega\times [0,T]; \mathcal{L}_2(K;H))$,
\[
I(f)\in L^p(\Omega; W^{\alpha, p}([0,T]; H))
\]
and there exists a constant $C(p,\alpha)>0$ independent of $f$ such that
\[
\mathbb{E}\|I(f)\|^p_{W^{\alpha, p}([0,T]; H)}\leq C(p, \alpha)\mathbb{E}\int^{T}_{0}\|f(t)\|^p_{\mathcal{L}_2(K;H)}dt.
\]
\end{lemma}
\begin{lemma}\label{lemma-3}
Let $B_0\subset B\subset B_1$ be Banach spaces, $B_0$ and $B_1$ reflexive, with compact embedding of $B_0$ in $B$. Let $p\in (1, \infty)$ and $\alpha \in (0, 1)$ be given. Let $X$ be the space
\[
X= L^p([0, T]; B_0)\cap W^{\alpha, p}([0,T]; B_1),
\]
endowed with the natural norm. Then the embedding of $X$ in $L^p([0,T];B)$ is compact.
\end{lemma}
\begin{lemma}\label{lemma-4}
If $B_1\subset \tilde{B}$ are two Banach spaces with compact embedding, and the real number $\alpha \in (0,1), p>1$ satisfy $\alpha p>1$, then the space $W^{\alpha, p}([0,T]; B_1)$ is compactly embedded into $C([0,T]; \tilde{B})$.
\end{lemma}

\begin{flushleft}
\textbf{Proof of Theorem \ref{thm-1}.}\quad The proof of the existence of weak solutions of (\ref{aa}) can be divided into three steps.

\textbf{Step 1.} \quad  Let $P_n$ be the operator from $D(A^{-\frac{3}{2}})$ to $D(A^{\frac{3}{2}})$ defined as
\[
P_n x=\sum^{n}_{i=1}\langle x, e_i\rangle e_i \quad x\in D(A^{-\frac{3}{2}}).
\]
Here, we denote by $\langle \cdot, \cdot \rangle$ the dual pairing between $D(A^{\frac{3}{2}})$ and $D(A^{-\frac{3}{2}})$.
Then
\[
\langle P_n x, y\rangle=\langle x, P_ny\rangle,
\]
for all $x,y \in D(A^{-\frac{3}{2}})$. Its restriction to $H$ is the orthogonal projection onto $P_n H:= \textrm{Span}\{e_1,\cdot\cdot\cdot, e_n\}$.
Let $B_n(Y,Y)$ be the Lipschitz operator in $P_nH$ defined as
\[
B_n(Y,Y)=\chi_n(Y)B(Y,Y) \quad Y\in P_nH,
\]
where $\chi_n: H\rightarrow \mathbb{R}$ is defined as $\chi_n(Y)=\Theta_n(|Y|)$,
with $\Theta_n: \mathbb{R}\rightarrow [0,1]$ of class $C^{\infty}$, such that
\[
\chi_n(Y)=\left\{
            \begin{array}{ll}
              1, & {\rm if} \ |Y|\leq n,\\
              0, & {\rm if} \  |Y|>n+1.
            \end{array}
          \right.
\]
Consider the classical Galerkin approximation scheme defined by
\begin{eqnarray}\label{aa-1}
\left\{
  \begin{array}{ll}
    dY_n+AY_ndt+P_nB_n(Y_n,Y_n)dt+P_nG(Y_n)dt=P_n\Psi(Y_n)dW(t), \quad t\in [0,T]\\
    Y_n(0)=P_ny_0.
  \end{array}
\right.
\end{eqnarray}
Noticing $B$ is locally Lipschitz from $V\times V$ to $D(A^{-\frac{3}{2}})$ and all the coefficients are continuous and linear growth in $P_n H$. Further, with the help of (\ref{a-2}) and (\ref{a-3}),  equation  (\ref{aa-1}) has a unique weak solution $Y_n\in L^2(\Omega; C([0,T];P_n H))$.
Applying It\^{o} formula to $|Y_n|^p$ for each $p\geq 2$,  by Lemma \ref{lemma-2}, we obtain that
 there exist two positive constants $C_1(p)$, $C_2$, which are independent of $n$, such that
\begin{eqnarray}\label{aa-2}
\mathbb{E}\left(\sup_{0\leq s\leq T}|Y_n(s)|^p\right)&\leq& C_1(p),\\
\label{aa-3}
\mathbb{E}\int^{T}_{0}\|Y_n(s)\|^2ds&\leq& C_2.
\end{eqnarray}

\textbf{Step 2.} \quad Decompose $Y_n$ as
\begin{eqnarray*}
Y_n(t)&=&P_ny_0-\int^t_0AY_n(s)ds-\int^t_0P_nB_n(Y_n(s),Y_n(s))ds
\\& &-\int^t_0P_nG(Y_n(s))ds+\int^t_0P_n\Psi(Y_n(s))dW(s)
\\  &:=& J^1_n+J^2_n(t)+J^3_n(t)+J^4_n(t)+J^5_n(t).
\end{eqnarray*}
Clearly, $\mathbb{E}|J^1_n|^2\leq C_3$. Since $\|AY_n\|_{V'}\leq \|Y_n\|$, we have
\begin{eqnarray*}
\mathbb{E}\|J^2_n(t)-J^2_n(s)\|^2_{V'}&\leq& C(t-s)\mathbb{E}\int^t_s\|AY_n(l)\|^2_{V'}dl\\
&\leq& C(t-s)\mathbb{E}\int^t_s\|Y_n(l)\|^2dl,
\end{eqnarray*}
which implies that
\[
\mathbb{E}\|J^2_n\|^2_{W^{1,2}([0,T]; V')}\leq C \mathbb{E}\int^{T}_{0}\|Y_n(s)\|^2ds\leq C_4.
\]
Refer to \cite{L-T-W-2}, it gives
\[
\|B(Y,Y_1)\|_{-3}\leq C|Y|\|Y_1\|,
\]
then
\begin{eqnarray*}
\mathbb{E}\|J^3_n(t)-J^3_n(s)\|^2_{D(A^{-\frac{3}{2}})}&\leq &C(t-s)\mathbb{E}\int^t_s\|B_n(Y_n,Y_n)\|^2_{-3}dl\\
&\leq & C(t-s)\mathbb{E}\sup_{t\in[0,T]}|Y_n(t)|^2\int^T_0\|Y_n(t)\|^2dt.
\end{eqnarray*}
Then, we deduce from (\ref{aa-2}) and (\ref{aa-3}) that
\[
\mathbb{E}\|J^3_n\|_{W^{1,2}([0,T]; D(A^{-\frac{3}{2}}))}\leq C_5\sqrt{C_1(2)C_2}.
\]
Notice that $\|G(Y_n)\|_{V'}\leq |Y_n|$, we obtain
\begin{eqnarray*}
\mathbb{E}\|J^4_n(t)-J^4_n(s)\|^2_{V'}&\leq &C(t-s)\mathbb{E}\int^t_s\|G(Y_n)\|^2_{V'}dl\\
&\leq & C(t-s)^2\mathbb{E}\sup_{t\in[0,T]}|Y_n(t)|^2.
\end{eqnarray*}
Hence, by (\ref{a-2}), we get
\[
\mathbb{E}\|J^4_n\|^2_{W^{1,2}([0,T]; V')}\leq C\mathbb{E}\left(\sup_{0\leq s\leq T}|Y_n(s)|^2\right)\leq C_6.
\]
Moreover, using Lemma \ref{lemma-5}, \textbf{Hypothesis H0} and (\ref{a-2}), we get
\[
\mathbb{E}\|J^5_n\|^2_{W^{\beta,2}([0,T]; H)}\leq C\mathbb{E}\left(\int^T_0\|P_n \psi(Y_n(s))\|^2_{\mathcal{L}_2(H;H)}ds\right)
\leq CT\mathbb{E}\sup_{0\leq s\leq T}|Y_n(s)|^2\leq C_7(\beta),
\]
for all $\beta\in (0,\frac{1}{2})$.
Since  $W^{1,p}([0,T]; B)\subset W^{\alpha,p}([0,T]; B)$
for all Banach space $B$, provided $\alpha \in (0,1)$ and $p>1$. Hence, combining all the previous inequalities, we obtain
\[
\mathbb{E}\|Y_n\|_{W^{\beta,2}([0,T]; D(A^{-\frac{3}{2}}))}\leq C_8(\beta),
\]
for all $\beta\in (0,\frac{1}{2})$. Recall (\ref{aa-3}), we know the laws $\mathcal{D}(Y_n)$ are bounded uniformly in probability in
\[
L^2([0,T];V)\cap W^{\beta,2}\left([0,T];D(A^{-\frac{3}{2}})\right).
\]
Thus, by Lemma \ref{lemma-3}, we conclude that the family $\mathcal{D}(Y_n)$ is tight in $L^2([0,T]; H)$.

Applying Lemma \ref{lemma-4} and by (\ref{aa-2}), similar to the term $J^5_n$, we have $\mathcal{D}(Y_n)$ is tight in $C\left([0,T]; D(A^{-\frac{\gamma}{2}})\right)$, for all given $\gamma>3$. Thus, we can find a subsequence still denoted by $Y_n$, such  that $\mathcal{D}(Y_n)$ converges weakly in $L^2([0,T];H)\cap C\left([0,T]; D(A^{-\frac{\gamma}{2}})\right)$.

\textbf{Step 3.} \quad  Fix $\gamma>3$. By Skorohod embedding theorem, there exists a stochastic basis  $(\Omega^1, \mathcal{F}^1,\{\mathcal{F}^1_t\}_{t\in [0,T]}, P^1)$ and $L^2([0,T];H)\cap C\left([0,T]; D(A^{-\frac{\gamma}{2}})\right)-$valued random variables $Y^1, Y^1_n, n\geq1$ on this basis, such that $Y^1_n$ has the same law of $Y_n$ on $L^2([0,T];H)\cap C\left([0,T]; D(A^{-\frac{\gamma}{2}})\right)$, and $Y^1_n\rightarrow Y^1$ in $L^2([0,T];H)\cap C\left([0,T]; D(A^{-\frac{\gamma}{2}})\right)$, $P^1-$a.s..  Of course, for each $n$,
\[
\mathcal{D}(Y^1_n)(C([0,T]; P_n H))=1,
\]
and by (\ref{aa-2}) and (\ref{aa-3}), we have
\begin{eqnarray*}
\mathbb{E}^{P^1}\left(\sup_{0\leq s\leq T}|Y^1_n(s)|^p\right)&\leq& C_1(p),\\
\mathbb{E}^{P^1}\int^T_0\|Y^1_n(s)\|^2ds&\leq& C_2,
\end{eqnarray*}
for all $n$ and  $p\geq 2$. Hence, we also have
\begin{equation}\notag
Y^1\in L^2([0,T]; V)\cap L^{\infty}([0, T]; H) \quad  P^1-a.s.
\end{equation}
and $Y^1_n\rightarrow Y^1$ weakly in $L^2(\Omega\times [0,T]; V)$.

For each $n\geq 1$, the process $M^1_n(t)$ with trajectories in $C([0,T];H)$ defined as
\[
M^1_n(t)=Y^1_n(t)-P_nY^1(0)+\int^t_0AY^1_n(s)ds+\int^t_0P_nB_n(Y^1_n(s),Y^1_n(s))ds+\int^t_0P_nG(Y^1_n(s))ds.
\]
In fact, $M^1_n$ is a square integrable martingale with respect to the filtration
\[
\Big\{\mathcal{G}^1_n\Big\}_t=\sigma\Big\{Y^1_n(s), s\leq t\Big\},
\]
with quadratic variation
\[
[M^1_n]_t=\int^t_0P_n\Psi(Y^1_n)\Psi(Y^1_n)^*P_nds.
\]
Then by a standard method (see \cite{F-G}), we obtain the existence of weak
solutions.
\end{flushleft}
$\hfill\blacksquare$

\begin{remark}
Refer to \cite{D-G-T}, the existence of the strong solution of 3D primitive equations (\ref{aa}) with initial data in $V$ is proved by Galerkin approximation method. Since $V\subset H$, the strong solution is also a weak solution which is a limitation of Galerkin approximations.
\end{remark}

%
\section{Exponential mixing for Galerkin sequences }\label{se-1}
In this section, we devote to proving exponential mixing for Galerkin sequences by coupling method. Firstly, we introduce some preliminaries.
\subsection{Preliminary}
Let $(\lambda_1,\lambda_2)$ be two distributions on a Polish space $(E, \mathcal{B}(E))$.  $(\Omega, \mathcal{F}, \mathbb{P})$ is a probability space. Denote by $(Z_1,Z_2)$  two random variables $(\Omega, \mathcal{F})\rightarrow (E, \mathcal{B}(E))$. We say that $(Z_1,Z_2)$ is a coupling of $(\lambda_1,\lambda_2)$ if $\lambda_{i}= \mathcal{D}(Z_i)$ for $i= 1, 2$. The total variation $\|\lambda\|_{var}$ of a finite real measure $\lambda$ on $E$ is defined as
 \[
 \|\lambda\|_{var}=\sup\Big\{|\lambda(\Gamma)| \ |\  \Gamma \in \mathcal{B}(E)\Big\},
\]
where  $\mathcal{B}(E)$ stands for the set of the Borelian subsets of $E$.

 The next result is fundamental in the coupling methods (the proof can be found in \cite{O-C}).
\begin{lemma}\label{le-1}
Let  $(\lambda_1,\lambda_2)$ be two probability measures on $(E,\mathcal{B}(E))$. Then
\[
\|\lambda_1- \lambda_2\|_{var}= \min \mathbb{P}(Z_1\neq Z_2).
\]
The minimum is taken over all couplings $(Z_1,Z_2)$ of  $(\lambda_1,\lambda_2)$. There exists a coupling which reaches the minimum value. It is called a maximal coupling.
\end{lemma}

For $N\in \mathbb{N}$, we denote by $P_N$ the eigenprojector of $A$ associated to the first $N$ eigenvalues. Let $(\Omega, \mathcal{F}, \mathbb{P})$ be a probability space and $W$ be a cylindrical Wiener process on $H$ for $\mathbb{P}$. Consider the following finite dimensional approximation of (\ref{aa}):
\begin{eqnarray}\label{aa-4}
\left\{
  \begin{array}{ll}
    dY_N+AY_Ndt+P_NB_N(Y_N,Y_N)dt+P_NG(Y_N)dt=P_N\Psi(Y_N)dW(t), \\
    P_NY(0)=P_Ny_0.
  \end{array}
\right.
\end{eqnarray}
From Theorem \ref{thm-1}, we know for any $y_0\in H$, (\ref{aa-4}) has a unique solution $Y_N=Y_N(\cdot, y_0)=(v_N(\cdot, y_0),S_N(\cdot, y_0))$.
Define
\[
(\mathcal{P}^N_t\psi)(y_0)=\mathbb{E}[\psi(Y_N(t,y_0))], \quad \rm for\ \psi\in {B}_b(P_N H).
\]
Since $Y_N(\cdot, y_0)$ satisfies the strong Markov property, we deduce that $(\mathcal{P}^N_t)_{t\in \mathbb{R}^{+}}$ is a Markov transition semigroup on $P_N H$.

\textbf{In the following, $P_N$ is omitted for simplicity.}

Applying It\^{o} formula to $|Y_N(\cdot, y_0)|^2$, it gives
\begin{align}\notag
d|Y_N|^2+2\|Y_N\|^2dt&=-2\Big(Y_N, B_N(Y_N,Y_N)\Big)dt-2\Big(Y_N,G(Y_N)\Big)dt
\\ \notag
& \quad  +2\Big(Y_N, \Psi(Y_N)dW\Big)+\|\Psi(Y_N)\|^2_{\mathcal{L}_2(H;H)}dt.
\end{align}
By integration by parts, we have
\begin{align}\notag
\Big(Y_N, B_N(Y_N,Y_N)\Big)=0,
\end{align}
With the aid of Lemma \ref{lemma-2}, we get
\begin{equation}\notag
\begin{aligned}
&d|Y_N|^2+2\|Y_N\|^2dt\\
&=2\Big(Y_N, \Psi(Y_N)dW\Big)+2\left(v_N,\frac{1}{\sqrt{C_0}}\int^z_{-1}\nabla_H S_N dz'\right)dt+\|\Psi(Y_N)\|^2_{\mathcal{L}_2(H;H)}dt\\
&\leq2\Big(Y_N, \Psi(Y_N)dW\Big)+\frac{2}{\sqrt{C_0}}|S_N|\|v_N\|dt+\|\Psi(Y_N)\|^2_{\mathcal{L}_2(H;H)}dt.
\end{aligned}
\end{equation}
Then, by $ab\leq \frac{1}{2}(a^2+b^2)$ and $\|Y_N\|^2=\|v_N\|^2+\|S_N\|^2$,  we have
\begin{equation}\notag
d|Y_N|^2+\|Y_N\|^2dt+\|S_N\|^2dt\leq 2(Y_N, \Psi(Y_N)dW)+\frac{8}{C_0}|S_N|^2dt+\kappa dt.
\end{equation}
Since $C_0\geq \frac{8}{\lambda_1}$, $\|S_N\|^2\geq \lambda_1|S_N|^2$, then $\frac{8}{C_0}|S_N|^2\leq \|S_N\|^2$. Thus
\begin{equation}\label{eq-16}
d|Y_N|^2+\|Y_N\|^2dt\leq 2(Y_N, \Psi(Y_N)dW)+\kappa dt.
\end{equation}
 By $\|Y_N\|^2 \geq \mu_1|Y_N|^2$, integrating ${\rm{e}}^{\mu_1 t}$ on both sides of (\ref{eq-16}) and taking expectation, we obtain
\begin{equation}\label{eq-17}
\mathbb{E}\left(|Y_N(t)|^2\right)\leq {\rm{e}}^{-\mu_1t}|y_0|^2+\frac{\kappa}{\mu_1}.
\end{equation}
Hence, applying the Krylov-Bogoliubov Criterion (see \cite{D-Z-2}), we obtain that $(\mathcal{P}^N_t)_t$ admits an invariant measure $\mu_N$ and that every invariant measure has a moment of order two in $H$.
 Let $Y^N_0$ be a random variable whose law is $\mu_N$ and which is independent of $W$, then $Y_N=Y_N(\cdot,Y^N_0)$ is a stationary solution of (\ref{aa-4}). Integrating (\ref{eq-16}), we obtain
\begin{equation}\label{eq-42}
\mathbb{E}|Y_N(t)|^2+\mathbb{E}\int^{t}_{0}\|Y_N(s)\|^2ds\leq\mathbb{E}|Y_N(0)|^2+\kappa t.
\end{equation}
Since the law of $Y_N(s)$ is $\mu_N$ for any $s\geq 0$, it follows that
\begin{equation}\label{aa-5}
\int_{P_N H}\|y\|^2\mu_N(dy)\leq \kappa.
\end{equation}

Similar to Theorem \ref{thm-1}, the laws $(\mathbb{P}^N_{\mu_N})$ of $Y_N(\cdot, Y^N_0)$ are tight in $L^2([0,T];H)\cap C([0,T]; D(A^{-\frac{\gamma}{2}}))$, $\gamma >3$. Then, for a subsequence, still denote by $(\mathbb{P}^N_{\mu_N})$, which converges in law to $\mathbb{P}_{\mu}$ a stationary solution of (\ref{aa}) with initial law $\mu$. We deduce from (\ref{aa-5}) that
\[
\int_{H}\|y\|^2\mu(dy)\leq\kappa.
\]

In general, we don't know whether $\mu$ is an invariant measure due to the lack of uniqueness and also we don't know whether the solution of (\ref{aa}) defines a Markov evolution. However, the above information is useful to prove the uniqueness of invariant measures.

\subsection{Exponential mixing for Galerkin sequences with initial value in $H$}
In the following, we aim to construct a coupling of solutions from $H$ satisfying exponential mixing property.

Consider solutions of (\ref{aa-4}) with initial data in $H$. Assume that \textbf{Hypothesis H1} holds. Let $N \in \mathbb{N}$ and $(y^1_0, y^2_0)\in H\times H$. Similar to \cite{K-S}, \cite{M}, we obtain that there exists a function $p_{N}(\cdot)>0$  such that
\begin{equation}\label{eq-22}
\|(\mathcal{P}^{N}_{1})^{*}\delta_{y^2_0}-(\mathcal{P}^{N}_{1})^{*}\delta_{y^1_0}\|_{var}\leq 1- p_N (|y^1_0|+|y^2_0|),
\end{equation}
where $\delta_{y^1_0}$ and $\delta_{y^2_0}$ are two Dirac measures on single point $y_0^1$ and $y_0^2$, respectively.
Applying Lemma \ref{le-1}, we construct a maximal coupling $(Z_1,Z_2)=\big(Z_1(y^1_0, y^2_0),Z_2(y^1_0, y^2_0)\big)$ of $\left((\mathcal{P}^{N}_{1})^{*}\delta_{y^1_0}, (\mathcal{P}^{N}_{1})^{*}\delta_{y^2_0}\right)$. It follows that
\begin{equation}\label{eq-18}
\mathbb{P}(Z_1=Z_2)\geq p_N (|y^1_0|+|y^2_0|)>0.
\end{equation}

Let $(W, \tilde{W})$ be a couple of independent cylindrical Wiener processes on $H\times H$ and $\delta >0$. Denote by $Y_N(\cdot, y^1_0)$ and $\tilde{Y}_{N}(\cdot, y^2_0)$ the solutions of (\ref{aa-4}) with initial data $y^1_0$ and $y^2_0$ associated to $W$ and $\tilde{W}$, respectively. Now we can construct a couple of random variables $(V_1, V_2)=(V_1(y^1_0, y^2_0), V_2(y^1_0, y^2_0))$ on $P_N H$ for $(Y_N, \tilde{Y}_N)$ as follows
\begin{eqnarray}\label{eq-19}
(V_1, V_2)=\left\{
             \begin{array}{ll}
              (Y_N(\cdot, y_0), Y_N(\cdot, y_0)), & {\rm if}\  y^1_0=y^2_0=y_0,
              \\

              (Z_1(y^1_0, y^2_0), Z_2(y^1_0, y^2_0)), & {\rm if}\  (y^1_0,y^2_0)\in B_{H\times H}(0, \delta)\backslash \{y^1_0=y^2_0\},
               \\

              (Y_N(\cdot, y^1_0), \tilde{Y}_N(\cdot, y^2_0)), & {\rm else},
            \end{array}
           \right.
\end{eqnarray}
where $B_{H\times H}(0, \delta)$ is the ball of $H\times H$ with radius $\delta$ . Then $(V_1(y^1_0, y^2_0), V_2(y^1_0, y^2_0))$ is a coupling of $\left((\mathcal{P}^{N}_{1})^{*}\delta_{y^1_0}, (\mathcal{P}^{N}_{1})^{*}\delta_{y^2_0}\right)$. It can be shown that it depends on $(y^1_0, y^2_0)$.
We then construct a coupling $(Y^1, Y^2)$ of $\left(\mathcal{D}(Y_N(\cdot, y^1_0)),\mathcal{D}(Y_N(\cdot, y^2_0))\right)$ by induction on $\mathbb{N}$. Firstly, setting $Y^{i}(0)=y^i_0$ for $i=1,2$, then assuming that we have constructed $(Y^1, Y^2)$ on $\{0,1,\cdot\cdot\cdot,k\}$. We take $(V_1,V_2)$ as above independent of $(Y^1,Y^2)$ and set
\[
Y^i(k+1)=V_i(Y^1(k), Y^2(k)) \quad \rm for\  i=1,2.
\]

Taking into account (\ref{eq-17}), it is easily shown that the time of return of $(Y^1, Y^2)$ in $B_{H\times H}(0,\frac{4\kappa}{\mu_1})$
admits an exponential moment. We choose $\delta= \frac{4\kappa}{\mu_1}$. It follows from (\ref{eq-18}) and (\ref{eq-19}) that $(Y^1(n), Y^2(n))\in B_{H\times H}(0,\delta)$ implies that the probability of $(Y^1, Y^2)$ having  coupled  at time $n+1$ is bounded below by $p_N(2\delta)>0$.
Finally, remark that if $(Y^1, Y^2)$ are coupled at time $n+1$, then they remain coupled for any time after. Combining these properties and using the fact that $(Y^1(n), Y^2(n))_{n\in \mathbb{N}}$ is a discrete strong Markov process, it is easily shown that
\begin{eqnarray}\label{eq-20}
\mathbb{P}(Y^1(n)\neq Y^2(n))\leq C_N {\rm{e}}^{-\gamma_{N}n}(1+|y^1_0|^2+|y^2_0|^2),
\end{eqnarray}
with $\gamma_{N}>0$.
Recall that $(Y^1, Y^2)$ is a coupling of $\left(\mathcal{D}(Y_N(\cdot, y^1_0)),\mathcal{D}(Y_N(\cdot, y^2_0))\right)$ on $\mathbb{N}$. It follows that $(Y^1(n), Y^2(n))$ is a coupling of $\left((\mathcal{P}^{N}_{n})^{*}\delta_{y^1_0}, (\mathcal{P}^{N}_{n})^{*}\delta_{y^2_0}\right)$. Combining Lemma \ref{le-1} and (\ref{eq-20}), we obtain, for $n\in \mathbb{N}$,
\[
\|(\mathcal{P}^{N}_{n})^{*}\delta_{y^2_0}- (\mathcal{P}^{N}_{n})^{*}\delta_{y^1_0}\|_{var}\leq C_N {\rm{e}}^{-\gamma_{N}n}(1+|y^1_0|^2+|y^2_0|^2).
\]

Setting $n=\lfloor t\rfloor$(the integer part of $t$) and integrating $(y^2_0, y^1_0)$ over $((\mathcal{P}^N_{t-n})^{*}\lambda)\bigotimes \mu_N$, where $\mu_N$ is an invariant measure, it follows that, for any $\lambda \in \mathcal{P}(P_N H)$ with $\int_{P_N H}|y|^2 \lambda(dy)<\infty$,
\begin{eqnarray}\label{eq-21}
\|(\mathcal{P}^N_{t})^{*}\lambda- \mu_{N}\|_{var}\leq C_N {\rm{e}}^{-\gamma_N t} \left(1+\int_{P_N H}|y|^2 \lambda(dy)\right).
\end{eqnarray}

This result (\ref{eq-21}) is useless when we consider (\ref{aa}), since the constants $C_N$, $\gamma_{N}$ in  (\ref{eq-21}) strongly depend on $N$.

 In the following, we devote to establishing a priori estimates uniformly in $N$ to ensure the exponential mixing property holds independent of $N$. Specifically, we obtain a uniform lower boundary for a coupling from $\mathbb{H}_3$ (see (\ref{eq-26}) in Proposition \ref{prp1} ), which is analogous to (\ref{eq-18}) but uniformly in $N$.



\subsection{Exponential mixing for Galerkin sequences with initial value in $\mathbb{H}_3$}\label{section-2}
As stated above, the initial data from $H$ is not enough to ensure some estimates uniformly in $N$. Here, we consider a coupling with initial data from $\mathbb{H}_3$.
Firstly, we obtain the following uniform lower bound.
\begin{prp}\label{prp1}
Assume that  \textbf{Hypothesis H1} holds. Then there exist $(\Upsilon,\delta)\in (0,1)\times (0,1)$ such that, for any $N \in \mathbb{N}$, there exists a coupling $(Z_1(y^1_0,y^2_0),Z_2(y^1_0,y^2_0))$ of $\left((\mathcal{P}^{N}_{\Upsilon})^{*}\delta_{y^1_0}, (\mathcal{P}^{N}_{\Upsilon})^{*}\delta_{y^2_0}\right)$ which measurably depends on $(y^1_0,y^2_0)\in \mathbb{H}_3\times \mathbb{H}_3$ and verifies
\begin{equation}\label{eq-26}
\mathbb{P}\left(Z_1(y^1_0,y^2_0)=Z_2(y^1_0,y^2_0)\right)\geq \frac{3}{4},
\end{equation}
provided
\begin{equation}\label{eq-27}
\|y^1_0\|^2_3 \vee \|y^2_0\|^2_3 \leq \delta.
\end{equation}
\end{prp}
The proof of Proposition \ref{prp1} is postponed to Sect. 6.

Assume that  \textbf{Hypothesis H1} holds. Let $N\in \mathbb{N}$ and $\Upsilon, \delta, Z_1, Z_2$ be the same as Proposition \ref{prp1}. Let $(W,\tilde{W})$ be a couple of independent cylindrical Wiener processes on $H\times H$. Denote by $Y_N(\cdot, y_0)$ and $\tilde{Y}_N(\cdot, y_0)$ the solution of (\ref{aa-1}) associated to $W$ and $\tilde{W}$, respectively. We construct a couple of random variables $(V_1,V_2)=(V_1(y^1_0,y^2_0),V_2(y^1_0,y^2_0))$ on $P_N H\times P_N H$ as follows
\begin{equation}\label{eq-61}
(V_1, V_2)=\left\{
             \begin{array}{ll}
              (Y_N(\cdot, y_0), Y_N(\cdot, y_0)), & {\rm if} \  y^1_0=y^2_0= y_0,  \\
              (Z_1(y^1_0, y^2_0), Z_2(y^1_0, y^2_0)), & {\rm if} \  (y^1_0,y^2_0)\in B_{\mathbb{H}_3\times \mathbb{H}_3}(0, \delta)\backslash \{y^1_0=y^2_0\},  \\
              (Y_N(\cdot, y^1_0), \tilde{Y}_N(\cdot, y^2_0)), & {\rm else}.
            \end{array}
           \right.
\end{equation}
 Then, we can construct $(Y^1, Y^2)$ by induction on $\Upsilon\mathbb{N}$. Indeed, firstly setting $Y^{i}(0)=y^{i}_0$ for $i=1,2.$ Then, assuming that we have constructed $(Y^1, Y^2)$ on $\{0,\Upsilon, 2\Upsilon,...,n\Upsilon\}$, taking $(V_1, V_2)$ as above independent of $(Y^1, Y^2)$ and setting
\begin{eqnarray}\label{equat-1}
Y^{i}((n+1)\Upsilon)= V_{i}(Y^1(n\Upsilon),Y^2(n\Upsilon)) \quad  {\rm for} \  i=1,2.
\end{eqnarray}
It follows that $(Y^1, Y^2)$ is a discrete strong Markov process and a coupling of $\big(\mathcal{D}(Y_N(\cdot, y^1_0)),\mathcal{D}(Y_N(\cdot, y^2_0))\big)$ on $\Upsilon\mathbb{N}$. Moreover, if $(Y^1,Y^2)$ are coupled at time $n\Upsilon$, then they remain coupled for any time after.

Define
\begin{equation}\label{eq-59}
\tau=\inf\left\{t\in \Upsilon\mathbb{N}\backslash\{0\}\mid \|Y^1(t)\|^2_3\vee\|Y^2(t)\|^2_3\leq \delta\right\}.
\end{equation}
\begin{prp}\label{pr-2}
Assume that  \textbf{Hypothesis H1} holds. There exist  $\alpha=\alpha(\Psi, \mathbb{T}^3, \varepsilon_0,\delta)>0$ and $K''= K''(\Psi, \mathbb{T}^3,\varepsilon_0,\delta)$ such that for any $(y^1_0, y^2_0)\in H\times H$ ,
\[
\mathbb{E}({\rm{e}}^{\alpha \tau})\leq K''(1+|y^1_0|^2+|y^2_0|^2),
\]
where $\varepsilon_0$ is given in \textbf{Hypothesis H1}.
\end{prp}


The proof of Proposition \ref{pr-2} is postponed to Sect. 6.



 Based on Propositions \ref{prp1} and \ref{pr-2} we can obtain the following exponential mixing property for Galerkin approximations.
\begin{prp}\label{pr-3}
Assume that  \textbf{Hypothesis H1} holds. Then there exist $C=C(\Psi, \mathbb{T}^3)>0$ and $\gamma=\gamma(\Psi, \mathbb{T}^3)>0$ such that for any $N\in \mathbb{N}$, there exists a unique invariant measure $\mu_N$ for $(\mathcal{P}^N_t)_{t\in \mathbb{R}^{+}}$. Moreover, for any $\lambda\in {P}(P_N H)$ with $\int_{P_N H}|y|^2 \lambda(dy)<\infty$,
\begin{eqnarray}\label{eq-23}
\|(\mathcal{P}^N_{t})^{*}\lambda- \mu_{N}\|_{var}\leq C {\rm{e}}^{-\gamma t} \left(1+\int_{P_N H}|y|^2 \lambda(dy)\right),
\end{eqnarray}
$\|\cdot\|_{var}$ is the total variation norm associated to the space $\mathbb{H}_s$, for $s<-3$.
\end{prp}

\begin{proof}
Based on the previous preparations, we can achieve this proposition using the same argument as \cite{O-C}.
Given $(y^1_0, y^2_0)\in \mathbb{H}_3\times \mathbb{H}_3$, the process $(Y^1, Y^2)$ is defined by (\ref{equat-1}).
Let $\delta>0$, $\Upsilon\in (0,1)$ be the same as Proposition \ref{prp1} and $\tau$ is defined by (\ref{eq-59}), set
\[
\tau_1=\tau, \quad \tau_{k+1}=\inf\left\{t>\tau_k\mid \|Y^1(t)\|^2_3\vee\|Y^2(t)\|^2_3\leq \delta\right\}.
\]
We deduce from the strong Markov property of $(Y^1, Y^2)$ and Proposition \ref{pr-2} that
\[
\mathbb{E}({\rm{e}}^{\alpha \tau_{k+1}})\leq K'' \mathbb{E}\left({\rm{e}}^{\alpha \tau_k}(1+|Y^1(\tau_k)|^2+|Y^2(\tau_k)|^2)\right),
\]
which yields,
\[
\left\{
  \begin{aligned}
   &\mathbb{E}({\rm{e}}^{\alpha \tau_{k+1}})\leq cK''(1+2\delta)\mathbb{E}({\rm{e}}^{\alpha \tau_{k}}) ,
   \\
   &\mathbb{E}( {\rm{e}}^{\alpha \tau_{1}}) \leq K''(1+|y^1_0|^2+|y^2_0|^2).
  \end{aligned}
\right.
\]
It follows that there exists $K>0$ such that
\[
\mathbb{E}({\rm{e}}^{\alpha \tau_{k}})\leq K^k (1+|y^1_0|^2+|y^2_0|^2).
\]
Hence, applying the Jensen's inequality, we obtain, for any $\theta\in (0,1)$,
\begin{equation}\label{eq-60}
\mathbb{E}({\rm{e}}^{\theta \alpha \tau_{k}})\leq K^{\theta k} (1+|y^1_0|^2+|y^2_0|^2).
\end{equation}
Taking into account Proposition \ref{prp1} and  (\ref{eq-61}) that
\[
\mathbb{P}\left(Y^1(\Upsilon)\neq Y^2(\Upsilon)\right)\leq \frac{1}{4},
\]
provided $(y^1_0, y^2_0)$ is in the ball of $\mathbb{H}_3\times \mathbb{H}_3$ with radius $\delta$.
Define
\[
k_0=\inf\left\{k\in \mathbb{N}\mid Y^1(\tau_k +\Upsilon)= Y^2(\tau_k +\Upsilon)\right\}.
\]
By strong Markov property of $(Y^1,Y^2)$, we have
\begin{equation}\label{eq-62}
\mathbb{P}(k_0>n)\leq \left(\frac{1}{4}\right)^n,
\end{equation}
which implies  $k_0 < \infty$ almost surely.
Let $\theta\in (0,1)$, we deduce from Cauchy-Schwarz inequality that
\[
\mathbb{E}({\rm{e}}^{\frac{\theta}{2} \alpha \tau_{k_0}})=\sum^{\infty}_{n=1}\mathbb{E}\left({\rm{e}}^{\frac{\theta}{2} \alpha \tau_n}I_{k_0=n}\right)
\leq\sum^{\infty}_{n=1}\sqrt{\mathbb{P}(k_0\geq n)\mathbb{E}({\rm{e}}^{\theta \alpha \tau_n})}.
\]
Combining (\ref{eq-60}) and (\ref{eq-62}),  we deduce
\[
\mathbb{E}({\rm{e}}^{\frac{\theta}{2} \alpha \tau_{k_0}})\leq \left(\sum^{\infty}_{n=0}\left(\frac{K^{\theta}}{2}\right)^n\right)(1+|y^1_0|^2+|y^2_0|^2).
\]
Hence, choosing $\theta\in (0,1)$ sufficiently small, we obtain that there exists $\gamma >0$
 independent of $N\in \mathbb{N}$ such that
\begin{equation}
\mathbb{E}({\rm{e}}^{\gamma\tau_{k_0}})\leq 4(1+|y^1_0|^2+|y^2_0|^2).
\end{equation}

Recall that if $(Y^1,Y^2)$ are coupled at time $t\in \Upsilon\mathbb{N}$, then they remain coupled for any time after. Hence, $Y^1(t)=Y^2(t)$ for $t>\tau_{k_0}$. It follows that
\[
\mathbb{P}(Y^1(n\Upsilon)\neq Y^2(n\Upsilon))\leq 4{\rm{e}}^{-\gamma n\Upsilon}(1+|y^1_0|^2+|y^2_0|^2).
\]
Since $\left(Y^1(n\Upsilon), Y^2(n\Upsilon)\right)$ is a coupling of $\left((\mathcal{P}^N_{n\Upsilon})^{*}\delta_{y^1_0},(\mathcal{P}^N_{n\Upsilon})^{*}\delta_{y^2_0}\right)$, we deduce from Lemma \ref{le-1} that
\begin{equation}\label{eq-71}
\|(\mathcal{P}^N_{n\Upsilon})^{*}\delta_{y^1_0}-(\mathcal{P}^N_{n\Upsilon})^{*}\delta_{y^2_0}\|_{var}\leq 4{\rm{e}}^{-\gamma n\Upsilon}(1+|y^1_0|^2+|y^2_0|^2),
\end{equation}
for any $n\in \mathbb{N}$ and any $(y^1_0,y^2_0)\in \mathbb{H}_3\times \mathbb{H}_3$.
Recall that the existence of an invariant measure $\mu_N \in {P}(P_N H)$ can be justified by (\ref{eq-42}). Let $\lambda \in {P}(H)$ and $t\in \mathbb{R}^{+}$. We set $n=\lfloor\frac{t}{\Upsilon}\rfloor$ and $C=4{\rm{e}}^{\gamma \Upsilon}$. Integrating $(y^1_0,y^2_0)$ over $((\mathcal{P}^N_{t-n\Upsilon})^{*}\lambda)\otimes \mu_N$
in (\ref{eq-71}), we obtain
\begin{equation}\label{eq-1000}
\|(\mathcal{P}^N_{t})^{*}\lambda-\mu_N\|_{var}\leq C{\rm{e}}^{-\gamma t}\left(1+\int_{H}|y|^2\lambda(dy)\right),
\end{equation}
which implies (\ref{eq-24}).

\end{proof}

Now, we are in the position to prove Theorem \ref{th-1}.

\noindent\textbf{Proof of Theorem \ref{th-1}. }\quad
Let $\lambda \in P(H)$ and $\mathbb{E}_{\lambda}$ be an expectation under the initial distribution $\lambda$. Since $\|\cdot\|_{var}$ is the dual norm of $|\cdot|_{\infty}$, we have for any finite measure $\lambda'$ on $\mathbb{H}_s$ for $s<-3$,
\begin{eqnarray}\label{l-1}
\|\lambda'\|_{var}=\sup_{|g|_{\infty}\leq 1}|\int_{\mathbb{H}_s}g(x)\lambda'(dx)|,
\end{eqnarray}
where the supermum is taken over $g\in UC_b(\mathbb{H}_s)$ which verifies $|g|_{\infty}<1$. Hence,
(\ref{eq-23}) is equivalent to
\begin{eqnarray}\label{eq-24}
\left|\mathbb{E}_{\lambda}\Big(g(Y_N(t))\Big)-\int_{P_N H}g(y)\mu_{N}(dy)\right|\leq C{\rm{e}}^{-\gamma t}|g|_{\infty}\left(1+\int_{ H}|y|^2 \lambda(dy)\right),
\end{eqnarray}
for any $g\in UC_b(\mathbb{H}_s)$.

Assume (\ref{eq-24}) holds. From Theorem \ref{thm-1}, we know that for any given initial law $\lambda\in P(H)$, there exists a subsequence $\{N^{'}_k\}_k$ such that $Y^{N^{'}_k}_{\lambda}$ converges in distribution to a weak solution $Y_{\lambda}$ of (\ref{aa}) in $C([0,T];\mathbb{H}_s)$, for $s<-3$. Moreover, as discussed in the above of Section \ref{se-1}, $(\mathbb{P}^N_{\mu_N})$ are tight in $L^2([0,T];H)\cap C([0,T]; D(A^{-\frac{\gamma}{2}}))$, $\gamma >3$.
Hence, there exists a subsequence $\{N_k\}_k$ of $\{N^{'}_k\}_k$ such that $(\mathbb{P}^{N_k}_{\mu_{N_k}})_k$ converges in law to $\mathbb{P}_{\mu^{(\lambda, \{N_k\})}}$, which is a stationary solution of (\ref{aa}) with initial law $\mu^{(\lambda, \{N_k\})}$. $\mu^{(\lambda, \{N_k\})}$ stands for $\mu$ depends on $\lambda$ and the sequence $\{N_k\}$. Taking the subsequence $\{N_k\}_k$ in both sides of (\ref{eq-24}) and letting $k\rightarrow \infty$, we have
 \begin{eqnarray}\label{eq-31}
 \left|\mathbb{E}_{\lambda}\Big(g(Y(t))\Big)-\int_{H}g(y)\mu^{(\lambda, \{N_k\})}(dy)\right|\leq C{\rm{e}}^{-\gamma t}|g|_{\infty}\left(1+\int_{ H}|y|^2 \lambda(dy)\right).
 \end{eqnarray}
Thus, by (\ref{l-1}) and (\ref{eq-31}), we conclude that Theorem \ref{th-1} holds with $\mu=\mu^{(\lambda, \{N_k\})}$.
$\hfill\blacksquare$

\section{Uniqueness of stationary probability measure}\label{l-3}
Given the initial law $\lambda\in P(H)$ and  Galerkin approximation sequence $\{N_k\}$, let $\mu^{(\lambda, \{N_k\})}$ be the stationary probability measure obtained in Theorem \ref{th-1}. In this part, we devote to proving $\mu^{(\lambda, \{N_k\})}$ is actually independent of $\lambda$ and $\{N_k\}$.



 \begin{flushleft}
 \textbf{Proof of Corollary \ref{cor-1}.}\quad Given any
initial value $y\in V$, by \cite{D-G-T-Z}, there exists a unique strong solution of (\ref{aa}) denoted by $Y(t,y)$. Based on  Theorem \ref{th-1},  $Y(t,y)$ has an invariant measure $\mu^{(\delta_{y})}$, which does not
depend on $\{N_k\}$.
We claim that
\begin{eqnarray}\notag
\mu^{(\delta_{y})}\equiv\mu^V \quad \forall y\in V.
\end{eqnarray}
Indeed, let
$y_i\in V (i=1,2)$ are two different initial values,
 $\mu^{(\delta_{y_1})}$ and $\mu^{(\delta_{y_2})}$ are the corresponding
invariant measures. Notice that
\begin{eqnarray}\notag
&& \left|\int_Hg(y)\mu^{(\delta_{y_1})}(dy)-\int_{H}g(y)\mu^{(\delta_{y_2})}(dy)\right|\\ \notag
 &\leq& \left|\int_Hg(y)\mu^{(\delta_{y_1})}(dy)-\mathbb{E}g(Y^{y_1}_t)\right|\\ \notag
& &\ + \left|\mathbb{E}g(Y^{y_1}_t)-\mathbb{E}g(Y^{y_2}_t)\right|
+ \left|\mathbb{E}g(Y^{y_2}_t)-\int_{H}g(y)\mu^{(\delta_{y_2})}(dy)\right|\\
\label{equ-5}
&:=& I_1(t)+I_2(t)+I_3(t),
\end{eqnarray}
we have
\begin{eqnarray*}\notag
I_2=\left|\mathbb{E}g(Y^{y_1}_t)-\mathbb{E}g(Y^{y_2}_t)\right|
&\leq& \left|\mathbb{E}g(Y^{y_1}_t)-\mathbb{E}g(Y^{N, y_1}_t)\right|\\ \notag
&&\ +\left|\mathbb{E}g(Y^{N, y_1}_t)-\mathbb{E}g(Y^{N, y_2}_t)\right|
+\left|\mathbb{E}g(Y^{N, y_2}_t)-\mathbb{E}g(Y^{y_2}_t)\right|\\ \notag
&:=& I^{1,N}_2(t)+I^{2,N}_2(t)+I^{3,N}_2(t).
\end{eqnarray*}
Thus,
\begin{eqnarray}\label{equ-14}
\left|\int_Hg(y)\mu^{(\delta_{y_1})}(dy)-\int_{H}g(y)\mu^{(\delta_{y_2})}(dy)\right|
 \leq I_1(t)+I^{1,N}_2(t)+I^{2,N}_2(t)+I^{3,N}_2(t)+I_3(t).
\end{eqnarray}
By (\ref{eq-71}) and (\ref{eq-31}), for any $\varepsilon>0$, there exists $t_0>0$ such that for any $t>t_0$
$$
I_1(t)+I^{2,N}_2(t)+I_3(t)\leq \varepsilon\ \ {\rm uniformly\ for\ N}.
$$
Fix $t>t_0$, and let $N\rightarrow\infty$, we obtain $I^{1,N}_2(t)+I^{3,N}_2(t)\rightarrow 0$. Thus,
let $t, N\rightarrow \infty$ in (\ref{equ-14}) , we get
\begin{eqnarray*}
\left|\int_Hg(y)\mu^{(\delta_{y_1})}(dy)-\int_{H}g(y)\mu^{(\delta_{y_2})}(dy)\right|=0,
\end{eqnarray*}
which implies
\begin{eqnarray}\label{equ-6}
\mu^{(\delta_y)}=\mu^{V}\quad \forall \ y\in V.
\end{eqnarray}
Now, we are ready to show that all weak solutions which are limits of Galerkin approximation share the same stationary measure. Let $y_0\in H$, $y_1\in V$, and $\lambda=\delta_{y_0}$, consider
\begin{eqnarray}\notag
 && \left|\int_Hg(y)\mu^{(\delta_{y_0},\{N_k\})}(dy)-\int_{H}g(y)\mu^{V}(dy)\right|\\ \notag
 &\leq&  \left|\int_Hg(y)\mu^{(\delta_{y_0},\{N_k\})}(dy)-\int_{P_{N_k}H}g(y)\mu^{(\delta_{y_0})}_{N_k}(dy)\right|\\ \notag
 &&\ + \left|\int_{P_{N_k}H}g(y)\mu^{(\delta_{y_0})}_{N_k}(dy)-\mathbb{E}g(Y^{N_k,y_0}_t)\right|+\left|\mathbb{E}g(Y^{N_k,y_0}_t)-\mathbb{E}g(Y^{N_k,y_1}_t)\right| \\ \notag
 &&\ +\left|\mathbb{E}g(Y^{N_k,y_1}_t)-\mathbb{E}g(Y^{y_1}_t)\right| + \left|\mathbb{E}g(Y^{y_1}_t)-\int_Hg(y)\mu^V(dy)\right|\\
 \label{equ-13}
 &:=& K^{N_k}_1+K^{N_k}_2(t)+K^{N_k}_3(t)+K^{N_k}_4(t)+K_5(t).
 \end{eqnarray}
 By (\ref{eq-71}), (\ref{eq-1000}), (\ref{eq-31}) and (\ref{equ-6}), we have $K^{N_k}_2(t)+K^{N_k}_3(t)+K^{N_k}_4(t)+K_5(t)\rightarrow 0$ uniformly for $\{N_k\}$,  as $t\rightarrow\infty$.  Let $N_k\rightarrow\infty$, we obtain $K^{N_k}_1\rightarrow 0$.
 Thus, let $t, N_k \rightarrow\infty$ in (\ref{equ-13}), we deduce that
\begin{equation}\notag
\int_Hg(y)\mu^{(\delta_{y_0},\{N_k\})}(dy)=\int_{H}g(y)\mu^{V}(dy),
\end{equation}
which implies $\mu^{(\delta_{y_0},\{N_k\})}=\mu^V$. That is, $\mu^{(\delta_{y_0},\{N_k\})}$ is independent of $\delta_{y_0}$  and $\{N_k\}$.

\end{flushleft}
$\hfill\blacksquare$

\section{Proof of Propositions \ref{prp1}-\ref{pr-2}}

\subsection{Proof of Proposition \ref{prp1} }

Assume \textbf{Hypothesis H1} holds. Let $\Upsilon\in (0,1)$. We can construct $\Big(Z_1(y^1_0,y^2_0),Z_2(y^1_0,y^2_0)\Big)$ as the maximal coupling of $(\mathcal{P}^{*}_{\Upsilon}\delta_{y^1_0},\mathcal{P}^{*}_{\Upsilon}\delta_{y^2_0})$ using Lemma \ref{le-1}. Measurable dependence on $(y^1_0,y^2_0)$ follows from a slight extension of Lemma \ref{le-1} (see \cite{O-C-1}, Remark A.1). Recall $\kappa$ is the constant defined in \textbf{Hypothesis H1}.
In order to establish Proposition \ref{prp1}, it is sufficient to prove that there exists $c( \kappa, \mathbb{T}^3)$ independent of $\Upsilon\in (0,1)$ and  $N \in \mathbb{N}$ such that
\begin{equation}\label{eq-27-1}
\|(\mathcal{P}^{N}_{\Upsilon})^{*}\delta_{y^2_0}- (\mathcal{P}^{N}_{\Upsilon})^{*}\delta_{y^1_0}\|_{var}\leq c( \kappa, \mathbb{T}^3)\sqrt{\Upsilon},
\end{equation}
provided
\begin{equation}\label{eq-28}
\|y^1_0\|^2_3 \vee \|y^2_0\|^2_3 \leq \kappa \Upsilon^3.
\end{equation}
Then it suffices to choose $\Upsilon\leq 1/(4c( \kappa, \mathbb{T}^3))^2$ and $\delta= \kappa \Upsilon^3$.

As $\|\cdot\|_{var}$ is the dual norm of $|\cdot|_{\infty}$, (\ref{eq-27-1}) is equivalent to
\begin{equation}\label{eq-29}
\left|\mathbb{E}\left(g(Y_N(\Upsilon, y^2_0))-g(Y_N(\Upsilon, y^1_0))\right)\right|\leq 8|g|_{\infty}c(\kappa, \mathbb{T}^3)\sqrt{\Upsilon},
\end{equation}
for any $g\in UC_b(P_N H)$. Due to  $C^1_b (P_N H)\subset U C_b (P_N H)$ is dense, it suffices to prove  (\ref{eq-29}) holds for any $N \in \mathbb{N}, \Upsilon\in (0,1)$ and $g \in C^1_b (P_N H)$ provided (\ref{eq-28}) holds.
\subsubsection{Energy estimates.}
For any process $Y=(v,S)$, define the $\mathbb{H}_2-$energy of $Y$ at time $t$ by
\[
E^{\mathbb{H}_2}_{Y}(t):=\|Y(t)\|^2_2 +\int^{t}_{0}\|Y(s)\|^2_3 ds.
\]
\begin{lemma}\label{le-3}
Assume that  \textbf{Hypothesis H1} holds. There exist $K_0= K_0( \mathbb{T}^3)$ and $c=c( \mathbb{T}^3)$ such that for any $\Upsilon \leq 1$ and any $N\in \mathbb{N}$, we have
\[
\mathbb{P}\left(\sup_{(0,\Upsilon)}E^{\mathbb{H}_2}_{Y_N(\cdot, y_0)}(t)> K_0\right)\leq c\left(1+\frac{\kappa}{K_0}\right)\sqrt{\Upsilon},
\]
provided $\|y_0\|^2_2\leq \kappa \Upsilon$.
\end{lemma}
\begin{proof}
Denote
\[
Y_N=Y_N(\cdot, y_0), v_N=v_N(\cdot, y_0), S_N=S_N(\cdot, y_0).
\]
Applying It\^{o} formula to $\|Y_N\|^2_2$ and by (\ref{aa-4}), we have
\begin{equation}\notag
d\|Y_N\|^2_2+2\|Y_N\|^2_3dt=dM_{\mathbb{H}_2} + I_{\mathbb{H}_2}dt +J_{\mathbb{H}_2}dt+\|P_N\Psi(Y_N)\|^2_{\mathcal{L}_2(H; \mathbb{H}_2)}dt,
\end{equation}
let
\[
I_{\mathbb{H}_2}=I^v_{\mathbb{H}_2}+I^S_{\mathbb{H}_2},\quad J_{\mathbb{H}_2}=J^v_{\mathbb{H}_2}+J^S_{\mathbb{H}_2},
\]
where
\[
\left\{
  \begin{aligned}
    &I^v_{\mathbb{H}_2}=-2\left(A^2_1v_N, (v_N\cdot \nabla_H)v_N+\Phi(v_N)\frac{\partial v_N}{\partial z}\right), \quad J^v_{\mathbb{H}_2}=-2\left(A^2_1v_N,fk\times v_N+\nabla_H p_b-\frac{1}{\sqrt{C_0}}\int^{z}_{-1}\nabla_H S_N dz'\right),
    \\
    &I^S_{\mathbb{H}_2}=-2\left(A^2_2S_N, (v_N\cdot \nabla_H)S_N+\Phi(v_N)\frac{\partial S_N}{\partial z}\right), \quad J^S_{\mathbb{H}_2}=0,
    \\
   & M_{\mathbb{H}_2}=2\int^{t}_{0}\left(A^2 Y_N(s),\Psi(Y_N(s))dW(s)\right).
  \end{aligned}
\right.
\]
By  Lemma \ref{lemm-2}, the H\"{o}lder inequality and the Young's inequality, we have
\begin{eqnarray*}\notag
I^v_{\mathbb{H}_2}&\leq& c\| v_N\|_3|\nabla v_N|_3|\nabla v_N|_6+c\| v_N\|_3\| v_N\|_2|v_N|_{\infty}+c\| v_N\|_3\| v_N\|_2|\frac{\partial v_N}{\partial_z}|_{\infty}+c\| v_N\|_3|\nabla v_N|_3||\nabla \frac{\partial v_N}{\partial_z}|_6\\
&\leq &c\| v_N\|^{\frac{3}{2}}_3\|v_N\|^{\frac{3}{2}}_2\leq \frac{1}{4}\|v_N\|^2_3+c\|v_N\|^6_2.
\end{eqnarray*}
Similarly, we obtain
\begin{eqnarray*}\notag
J^v_{\mathbb{H}_2}&\leq&2\| v_N\|_3\| v_N\|+\frac{2}{\sqrt{C_0}}\| v_N\|_3\|S_N\|_2\leq \frac{1}{4}\|v_N\|^2_3+c\|v_N\|^2+c\|S_N\|^2_2,
\\
I^S_{\mathbb{H}_2}&\leq &\frac{1}{2}\|S_N\|^2_3+ c\|v_N\|^4_2+c\|v_N\|^8_2+c\|S_N\|^8_2+c\|S_N\|^4_2.
\end{eqnarray*}
Based on the above inequalities, we obtain
\begin{equation}\notag
d\|Y_N\|^2_2+\frac{3}{2}\|Y_N\|^2_3dt\leq c\|Y_N\|^8_2dt + c\kappa dt+ dM_{\mathbb{H}_2},
\end{equation}
where $\kappa$ is defined in \textbf{Hypothesis H1}. Then, we get
\begin{equation}\notag
d\|Y_N\|^2_2+\|Y_N\|^2_3dt\leq c\|Y_N\|^2_2\left(\|Y_N\|^6_2-8K^3_0\right)dt + c\kappa dt+ dM_{\mathbb{H}_2},
\end{equation}
for
\begin{equation}\label{eq-51}
K_0=\sqrt[3]{\frac{\mu_1}{16c}}.
\end{equation}
Set
\[
\sigma_{\mathbb{H}_2}=\inf\{t\in(0,\Upsilon)\mid \|Y_N(t)\|^2_2>2K_0\}.
\]
Since $\|y_0\|^2_2\leq \kappa \Upsilon$, we deduce that for any $t\in (0,\sigma_{\mathbb{H}_2})$,
\begin{equation}\label{eq-32}
\mathbb{E}^{\mathbb{H}_2}_{Y_N}(t)\leq M_{\mathbb{H}_2}(t)+c\kappa \Upsilon.
\end{equation}
From  \textbf{Hypothesis H1}, we know that $\Psi(y)^{*}A$ is bounded in $\mathcal{L}(\mathbb{H}_2;\mathbb{H}_2)$ by $c\kappa$. It follows that for any $t\in (0, \sigma_{\mathbb{H}_2})$,
\[
\langle M_{\mathbb{H}_2}\rangle(t)=4\int^{t}_{0}|P_N\Psi(Y_N(s))^{*}A^2Y_N(s)|^2ds \leq c\kappa \int^{t}_{0}\|Y_N\|^2_2ds\leq 2c\kappa K_0\Upsilon.
\]
Applying Burkholder-Davis-Gundy inequality, we have
\[
\mathbb{E}\left(\sup_{(0,\sigma_{\mathbb{H}_2})}M_{\mathbb{H}_2}(t)\right)\leq c\mathbb{E}\sqrt{\langle M_{\mathbb{H}_2}\rangle(\sigma_{\mathbb{H}_2})}\leq c\sqrt{K_0\kappa\Upsilon}\leq c(K_0 +\kappa)\sqrt{\Upsilon}.
\]
We deduce from (\ref{eq-32}) and $\Upsilon\leq 1$ that
\[
\mathbb{E}\left(\sup_{(0,\sigma_{\mathbb{H}_2})}E^{\mathbb{H}_2}_{Y_{N}}(t)\right)\leq c(K_0 +\kappa)\sqrt{\Upsilon},
\]
which yields
\[
\mathbb{P}\left(\sup_{(0,\sigma_{\mathbb{H}_2})}E^{\mathbb{H}_2}_{Y_{N}}(t)>K_0\right)\leq c(1+\frac{\kappa}{K_0})\sqrt{\Upsilon}.
\]

Let
$B=\left\{\sup_{(0,\sigma_{\mathbb{H}_2})}E^{\mathbb{H}_2}_{Y_{N}}(t)\leq
K_0\right\}$,
$A=\left\{\sup_{(0,\Upsilon)}E^{\mathbb{H}_2}_{Y_{N}}(t)\leq
K_0\right\}$. Since
$\sup_{(0,\sigma_{\mathbb{H}_2})}E^{\mathbb{H}_2}_{Y_{N}}(t)\leq
K_0$ implies $\sigma_{\mathbb{H}_2}=\Upsilon$, we have $B\subset A$, then $\mathbb{P}(A^c)\leq \mathbb{P}(B^c)$,
which implies Lemma \ref{le-3}.
\end{proof}
\subsubsection{Derivative estimates.}
Let $N\in \mathbb{N}$ and $(y_0, h)\in (\mathbb{H}_3)^2$, where $y_0= (v_0, S_0)$. Consider
\begin{equation}\label{eq-33}
\left\{
   \begin{aligned}
\frac{\partial \beta_N}{\partial t}+P_N(v_N\cdot \nabla_H)\beta_N&+P_N\Phi(v_N)\frac{\partial \beta_N}{\partial z}+P_N(\eta_N\cdot \nabla_H)Y_N+P_N\Phi(\eta_N)\frac{\partial Y_N}{\partial z}\\
&+
P_N\tilde{G}(\beta_N)                
+
A              
\beta_N            
=P_N\Psi'(Y_N)\beta_N
\frac{dW}{dt},
\\
\beta_N(s,s,y_0)\cdot h=P_N h,&              
\end{aligned}  \\
\
\right.
\end{equation}
where
\[
\tilde{G}(\beta_N)=\left(                 
  \begin{array}{c}   
  P^1_N fk\times\eta_N-\frac{1}{\sqrt{C_0}}P^1_N\int^{z}_{-1}\nabla_H \gamma_N dz'  \\  
  0  \\  
  \end{array}
\right) ,
\]
and
\[
\eta_N (t)=\eta_N(t,s,y_0)\cdot h, \quad \gamma_N (t)=\gamma_N(t,s,y_0)\cdot h \quad {\rm for}\  t\geq s.
\]
Denote $\beta_N=(\eta_N, \gamma_N)$ and $\beta_N (t)=\beta_N(t,s,y_0)\cdot h$.
The existence and uniqueness of the solutions of (\ref{eq-33}) is easily obtained.
Moreover, if $g\in C^1_b(P_N H)$, then, for any $t\geq 0$, we have
\begin{equation}\notag
\left(\nabla \left(\mathcal{P}^{N}_{t}g\right)(y_0), h\right)=\mathbb{E}\left(\nabla g(Y_N(t,y_0)),\beta_{N}(t,0,y_0)\cdot h\right).
\end{equation}
For process $Y=(v, S)$, set
\begin{equation}\label{eq-38}
\sigma(Y)=\inf\left\{t\in (0,\Upsilon) | \int^{t}_{0}\|Y(s)\|^2_3ds\geq K_0 +1\right\},
\end{equation}
where $K_0$ is defined by (\ref{eq-51}).
\begin{lemma}\label{le-4}
Assume that  \textbf{Hypothesis H1} holds. Then there exists $c=c(\kappa, \mathbb{T}^3)$ such that for any $N\in \mathbb{N}, \Upsilon\leq 1$ and $(y_0, h)\in (\mathbb{H}_3)^2$,
\[
\mathbb{E}\int^{\sigma(Y_N(\cdot, y_0))}_{0} \|\beta_N(t,0,y_0)\cdot h\|^2_4dt\leq c\|h\|^2_3.
\]
\end{lemma}

 \begin{proof} Set
 \begin{eqnarray*}
  \beta_N(t)&=&\beta_N (t,0, y_0)\cdot h,\ \eta_N(t)=\eta_N (t,0, y_0)\cdot h,\\
   \gamma_N(t)&=&\gamma_N (t,0, y_0)\cdot h\ , \sigma=\sigma(Y_N(\cdot, y_0)).
  \end{eqnarray*}
  Applying It\^{o} formula to $\|\beta_N(t)\|^2_3$ and by (\ref{eq-33}), it gives
\begin{equation}\label{eq-35}
d\|\beta_N(t)\|^2_3 +2 \|\beta_N(t)\|^2_4dt= dM_{\beta_N} +I_{\beta_N}dt+J_{\beta_N}dt +\|P_N(\Psi'(Y_N)\cdot\beta_N)\|^2_{\mathcal{L}_2(H;\mathbb{H}_3)}dt,
\end{equation}
where
\begin{equation}\notag
\left\{
        \begin{aligned}
         & M_{\beta_N}(t)=2\int^{t}_{0}\Big(A^3\beta_N, \big(P_N\Psi'(Y_N)\cdot\beta_N\big)dW\Big),
         \\
          &I_{\beta_N}=-2\left(A^3\beta_N,(v_N\cdot \nabla)\beta_N +(\eta_N \cdot \nabla)Y_N\right)-2\left(A^3\eta_N,\Phi(v_N)\frac{\partial {\beta_N}}{\partial z}+\Phi(\eta_N)\frac{\partial Y_N}{\partial z}\right),
          \\
          &J_{\beta_N}=-2\left(A^3\beta_N,\tilde{G}(\beta_N)\right).
        \end{aligned}
      \right.
\end{equation}
Let
\[
I_{\beta_N}=I_{\eta_N}+I_{\gamma_N},\quad J_{\beta_N}=J_{\eta_N}+J_{\gamma_N},
\]
where
\begin{equation}\notag
\left\{
        \begin{aligned}
          &I_{\eta_N}=-2\left(A^3_1\eta_N,(v_N\cdot \nabla_H)\eta_N+\Phi(v_N)\frac{\partial {\eta_N}}{\partial z} \right)-2\left(A^3_1\eta_N,(\eta_N \cdot \nabla_H)v_N+\Phi(\eta_N)\frac{\partial v_N}{\partial z}\right),
          \\
          &J_{\eta_N}=-2\left(A^3_1\eta_N,fk\times \eta_N-\frac{1}{\sqrt{C_0}}\int^{z}_{-1}\nabla_H \gamma_N dz'\right),
        \end{aligned}
      \right.
\end{equation}
and
\begin{equation}\notag
\left\{
        \begin{aligned}
          &I_{\gamma_N}=-2\left(A^3_2\gamma_N,(v_N\cdot \nabla_H)\gamma_N +\Phi(v_N)\frac{\partial {\gamma_N}}{\partial z}\right)-2\left(A^3_2\gamma_N,(\eta_N \cdot \nabla_H)S_N+\Phi(\eta_N)\frac{\partial S_N}{\partial z}\right),
          \\
          &J_{\gamma_N}=0.
        \end{aligned}
      \right.
\end{equation}
By Lemma \ref{lemm-2}, the H\"{o}lder inequalities, Sobolev embedding and the Young's inequality, we obtain
\begin{eqnarray*}
I_{\eta_N}&\leq& c\|\eta_N\|_4\|\eta_N\|_3|v_N|_{\infty}+c\|\eta_N\|_4|A_1 v_N|_3|\nabla_H \eta_N|_6
+c\|\eta_N\|_4\|v_N\|_3|\frac{\partial \eta_N}{\partial_z}|_{\infty}+c\|\eta_N\|_4|A_1\frac{\partial \eta_N}{\partial_z}||\nabla v_N|_{\infty}\\
&&\ +c\|\eta_N\|_4\|v_N\|_3|\eta_N|_{\infty}+c\|\eta_N\|_4|A_1 \eta_N|_3|\nabla_H v_N|_6
+c\|\eta_N\|_4\|\eta_N\|_3|\frac{\partial v_N}{\partial_z}|_{\infty}+c\|\eta_N\|_4|A_1\frac{\partial v_N}{\partial_z}||\nabla \eta_N|_{\infty}\\
&\leq & c\|\eta_N\|_4\|v_N\|_3\|\eta_N\|_3
\leq  \frac{1}{4}\|\eta_N\|^2_4+c\|v_N\|^2_3\|\eta_N\|^2_3.
\end{eqnarray*}
Using the H\"{o}lder inequalities, Sobolev embedding and the Young's inequality, we get
\begin{eqnarray*}
J_{\eta_N}&\leq & c\|\eta_N\|_4\|\eta_N\|_2+c\|\eta_N\|_4\|\gamma_N\|_3
\leq  \frac{1}{4}\|\eta_N\|^2_4+c\|\eta_N\|^2_2+c\|\gamma_N\|^2_3. \notag
\end{eqnarray*}
Similarly, we obtain
\[
I_{\gamma_N}+J_{\gamma_N}\leq \frac{1}{2}\|\gamma_N\|^2_4+c\|\gamma_N\|^2_3\|v_N\|^2_3+c\|\eta_N\|^2_3\|S_N\|^2_3.
\]

Collecting all the above inequalities, we get
\begin{equation}\label{eq-64}
d\|\beta_N\|^2_3+\frac{3}{2}\|\beta_N\|^2_4dt\leq c\|\beta_N\|^2_3(1+\|Y_N\|^2_3)dt+dM_{\beta_N}.
\end{equation}
Integrating and taking the expectation, we deduce that
\begin{equation}\notag
\mathbb{E}\left(\mathcal{E}(\sigma,0)\|\beta_N(\sigma)\|^2_3+\int^{\sigma}_{0}\mathcal{E}(t,0)\|\beta_N(t)\|^2_4dt\right)\leq \|h\|^2_3,
\end{equation}
where
\[
\mathcal{E}(t,0)={\rm{e}}^{-ct-c\int^{t}_{0}\|Y_N(r)\|^2_3dr}.
\]

From the definition of $\sigma$, we deduce that
\begin{equation}\notag
\mathbb{E}\int^{\sigma}_{0}\|\beta_N(t)\|^2_4dt \leq \|h\|^2_3 \exp\Big(c(K_0+1)+c\Upsilon \Big),
\end{equation}
which yields Lemma \ref{le-4}.
\end{proof}

Now, we are ready to prove Proposition \ref{prp1}.
\begin{flushleft}
\textbf{Proof of Proposition \ref{prp1}.} \quad As explained above, it suffices to prove (\ref{eq-29}) holds for any $N\in \mathbb{N}$, $\gamma\in (0,1)$ and $g\in C^1_b (P_N H)$ provided (\ref{eq-28}) holds.
Let $\psi\in C^{\infty}(\mathbb{R};[0,1])$ be defined by
\[
\psi=\left\{
       \begin{array}{ll}
         0, &\rm on \  (K_0+1, \infty), \\
         1, & \rm on \ (-\infty, K_0),
       \end{array}
     \right.
\]
where $K_0$ is defined by (\ref{eq-51}). For the process $Y$, set
\[
\psi_{Y}=\psi\left(\int^{\Upsilon}_{0}\|Y(s)\|^2_3ds\right).
\]
Notice that
\begin{equation}\label{eq-40}
\left|\mathbb{E}\left(g(Y_N(\Upsilon,y^2_0))-g(Y_N(\Upsilon,y^1_0))\right)\right|\leq I_0 + |g|_{\infty}(I_1+I_2),
\end{equation}
where
\[
\left\{
  \begin{aligned}
    &I_0=\left|\mathbb{E}\left(g(Y_N(\Upsilon ,y^2_0))\psi_{Y_N(\cdot, y^2_0)}-g(Y_N(\Upsilon,y^1_0))\psi_{Y_N(\cdot, y^1_0)}\right)\right|,
    \\
    &I_i=\mathbb{P}\left(\int^{\Upsilon}_{0}\|Y_N(s,y^i_0)\|^2_3 ds > K_0\right).
  \end{aligned}
\right.
\]
For $\theta \in [1,2]$, set
\[
\left\{
  \begin{aligned}
    &y^{\theta}_0=(2-\theta)y^1_0 + (\theta-1)y^2_0, \quad Y_{\theta}= Y_{N}(\cdot, y^{\theta}_0), \quad v_{\theta}= v_{N}(\cdot, y^{\theta}_0), \quad S_{\theta}= S_{N}(\cdot, y^{\theta}_0),
    \\
   &\beta_{\theta}(t)=\beta_{N}(t,0, y^{\theta}_0),\quad \eta_{\theta}(t)=\eta_{N}(t,0, y^{\theta}_0), \quad \sigma_{\theta}=\sigma(Y_{\theta}),
  \end{aligned}
\right.
\]
 where $\sigma$ is defined by (\ref{eq-38}). For better readability, the dependence on $N$ has been omitted.
Setting
\[
h=y^2_0-y^1_0,
\]
 we have
\begin{equation}\notag
I_0\leq \int^{2}_{1}|J_{\theta}|d \theta , \quad J_{\theta}=\left(\nabla_{y^{\theta}_0} \mathbb{E}(g(Y_{\theta}(\Upsilon))\psi_{Y_{\theta}}),h\right).
\end{equation}

To estimate $J_{\theta}$, applying a truncated Bismut-Elworthy-Li formula, similar to \cite{O-C}, we have
\begin{equation}\notag
J_{\theta}=\frac{1}{\Upsilon}J^{'}_{\theta, 1}+2J'_{\theta,2},
\end{equation}
 where
\begin{equation}\notag
\left\{
  \begin{aligned}
  & J^{'}_{\theta, 1}=\mathbb{E}\left[g\left(Y_{\theta}(\Upsilon)\right)\psi_{Y_{\theta}}\int^{\sigma_{\theta}\wedge \Upsilon}_{0}\left(\Psi^{-1}(Y_{\theta}(t))\cdot\beta_{\theta}(t)\cdot h, dW(t)\right)\right] ,
  \\
   & J^{'}_{\theta, 2}=\mathbb{E}\left[g(Y_{\theta}(\Upsilon))\psi^{'}_{Y_{\theta}}\int^{\sigma_{\theta}\wedge \Upsilon}_{0}(1-\frac{t}{\Upsilon})\left(A^{\frac{3}{2}} Y_{\theta}(t),A^{\frac{3}{2}}(\beta_{\theta}(t)\cdot h)dt\right)\right],
   \\
  &\psi^{'}_{Y_{\theta}}=\psi^{'}\left(\int^{\Upsilon}_{0}\|Y_{\theta}(s)\|^2_3ds\right).\\
\end{aligned}
\right.
\end{equation}
It follows from  \textbf{Hypothesis H1} that
\[
|J^{'}_{\theta, 1}|\leq |g|_{\infty}\kappa\sqrt{\mathbb{E}\int^{\sigma_{\theta}\wedge \Upsilon}_{0}\|\beta_{\theta}(t)\cdot h\|^2_4 dt},
\]
and from H\"{o}lder inequality that
\[
|J^{'}_{\theta, 2}|\leq |g|_{\infty}|\psi^{'}|_{\infty}\sqrt{\mathbb{E}\int^{\sigma_{\theta}\wedge \Upsilon}_{0}\|Y_{\theta}(t)\|^2_3dt}\sqrt{\mathbb{E}\int^{\sigma_{\theta}\wedge \Upsilon}_{0}\|\beta_{\theta}(t)\cdot h\|^2_3dt}.
\]
Hence for any $\Upsilon<1$,
\begin{equation}\label{eq-39}
|J_{\theta}|\leq c(\kappa, \mathbb{T}^3)|g|_{\infty}\frac{1}{\Upsilon}\sqrt{\mathbb{E}\int^{\sigma_{\theta}\wedge \Upsilon}_{0}\|\beta_{\theta}(t)\cdot h\|^2_4dt}.
\end{equation}

Combining (\ref{eq-39}) and Lemma \ref{le-4}, we have
\begin{equation}\notag
|J_{\theta}|\leq c(\kappa, \mathbb{T}^3)|g|_{\infty}\frac{\|h\|_3}{\Upsilon},
\end{equation}
which yields,
\[
I_0\leq c(\kappa, \mathbb{T}^3)|g|_{\infty}\sqrt{\Upsilon}.
\]
Since $\kappa \Upsilon^3\leq \kappa \Upsilon$, we can apply Lemma
\ref{le-3} to control $I_1+I_2$ in (\ref{eq-40}) if (\ref{eq-28})
holds. Hence (\ref{eq-29}) follows provided (\ref{eq-28}) holds,
which yields Proposition \ref{prp1}.
\end{flushleft}

$\hfill\blacksquare$

\subsection{Proof of Proposition \ref{pr-2}}
In order to prove Proposition \ref{pr-2}, we need to verify the following Lemmas \ref{le-5}-\ref{le-9}. Firstly, similar to \cite{O-C}, we have
\begin{lemma}\label{le-5}
Assume that  \textbf{Hypothesis H1} holds. For any $t, M>0$, there exists $p_0(t,M)=p_0(t,M,\varepsilon_0, \{|\Psi_n|_{\infty}\}_n, \mathbb{T}^3)>0$ such that for any adapted process $Y$,
\[
\mathbb{P}\left(\sup_{(0,t)}\|Z(s)\|^2_3\leq M\right)\geq p_0(t,M),
\]
where
\[
Z(t)=\int^{t}_{0}{\rm{e}}^{-A(t-s)}\Psi(Y(s))dW(s).
\]
\end{lemma}
Using this estimation, we can estimate the moment of the first time in a small ball in $H$. Let $\delta_3>0$, set
\[
\tau_{L^2}=\tau \wedge \inf \left\{t\in \Upsilon\mathbb{N}\mid |Y^1(t)|^2\vee|Y^2(t)|^2\leq \delta_3\right\}.
\]
\begin{lemma}\label{le-6}
Assume that  \textbf{Hypothesis H1} holds. Then, for any $\delta_3>0$, there exist $C_3(\delta_3)$ and $\gamma_3(\delta_3) $ such that for any $(y^1_0,y^2_0)\in H \times H$,
\[
\mathbb{E}({\rm{e}}^{\gamma_3 \tau_{L^2}})\leq C_3\left(1+|y^1_0|^2+|y^2_0|^2\right).
\]
\end{lemma}

\begin{proof}
Recall (\ref{eq-17}), we have
\[
\mathbb{E}|Y_N(t)|^2\leq {\rm{e}}^{-\mu_1t}|y_0|^2+\frac{\kappa}{\mu_1}.
\]
Since $(Y^1, Y^2)$ is a coupling of $\left(\mathcal{D}(Y_N(\cdot,y^1_0)),\mathcal{D}(Y_N(\cdot,y^2_0))\right)$ on $\Upsilon\mathbb{N}$, we obtain
\begin{equation}\notag
\mathbb{E}\left(|Y^1(n\Upsilon)|^2+|Y^2(n\Upsilon)|^2\right)\leq {\rm{e}}^{-\mu_1 n\Upsilon}(|y^1_0|^2+|y^2_0|^2)+2\frac{\kappa}{\mu_1}.
\end{equation}
Since $(Y^1,Y^2)$ is a strong Markov process, it can be deduced that there exist $C_7$ and $\gamma_7$ such that
\begin{equation}\label{eq-45}
\mathbb{E}\left({\rm{e}}^{\gamma_7 \tau'_{L^2}}\right)\leq C_7(1+|y^1_0|^2+|y^2_0|^2),
\end{equation}
where
\[
\tau'_{L^2}=\inf\left\{t\in \Upsilon\mathbb{N}\backslash \{0\}\mid |Y^1(t)|^2+|Y^2(t)|^2\leq 4\kappa\right\}.
\]
Taking (\ref{eq-45}) into account, in order to establish Lemma \ref{le-6}, it is sufficient to prove that there exist $(p_8(\delta_3,t),\Upsilon_8(\delta_3)) $ such that
\begin{equation}\label{eq-46}
\mathbb{P}\left(|Y_N(t,y_{0})|^2\leq \delta_3\right)\geq p_8>0,
\end{equation}
provided $N\in \mathbb{N}, t\geq \Upsilon_8(\delta_3)$, $|y_0|^2\leq 4\kappa$ and $\Upsilon_8(\delta_3)$ is independent of $y_0$.
Set
\[
Z(t)=\int^{t}_{0}{\rm{e}}^{-A(t-s)}\Psi(Y(s))dW(s), \quad X_N=Y_N-P_NZ,
\quad
 N(\omega)=\sup_{(0,t)}\|Z(s,\omega)\|^2_3 \quad {\rm for}\ \omega\in \Omega.
\]

Assume  that there exist $M_8(\delta_3)>0$ and $\Upsilon_8(\delta_3)$ such that for $\omega\in \Omega$,
\begin{equation}\label{eq-47}
N(\omega)\leq M_8(\delta_3)\wedge \frac{\delta_3}{4}\ {\rm implies} \quad |X_N(t,\omega)|^2\leq \frac{\delta_3}{4},
\end{equation}
provided $t\geq \Upsilon_8(\delta_3)$ and $|y_{0}|^2\leq 4\kappa$ . Then, we deduce (\ref{eq-46}) holds from Lemma \ref{le-5}   with $M= M_8(\delta_3)\wedge \frac{\delta_3}{4}$.

We now prove (\ref{eq-47}).
From (\ref{aa-4}), we know
\begin{eqnarray*}
\left\{
   \begin{aligned}
\frac{\partial X_N}{\partial t}&+P_N(v_N\cdot \nabla_H)(X_N+P_N Z)+P_N\Phi(v_N)\frac{\partial (X_N+P_N Z)}{\partial z}+
P_NG(X_N+P_N Z)                
+
A                 
X_N            
=0,
\\
X_N(0)&=P_N y_0.
\end{aligned}
\right.
\end{eqnarray*}
Let
 \[
 X_N=(\omega_N, g_N)=(v_N, S_N)-(P^1_NZ_1,P^2_NZ_2),
 \]
where
\[
Z_1(t)=\int^{t}_{0}{\rm{e}}^{-A_1(t-s)}\phi(v_N(s),S_N(s))dW_1(s),\quad Z_2(t)=\int^{t}_{0}{\rm{e}}^{-A_2(t-s)}\varphi(v_N(s),S_N(s))dW_2(s).
\]
For $\omega\in \Omega$, setting
\[
N_1(\omega)=\sup_{(0,t)}\|Z_1(s,\omega)\|^2_3, \quad N_2(\omega)=\sup_{(0,t)}\|Z_2(s,\omega)\|^2_3,
\]
we have
\[
N_1(\omega)\vee N_2(\omega)\leq N(\omega).
\]

From (\ref{aa-4}), we have
\begin{eqnarray}\notag
 \frac{\partial \omega_N}{\partial t}&+&\Big((\omega_N+ Z_1)\cdot \nabla_H\Big)(\omega_N+ Z_1)+\Phi(\omega_N+Z_1)\frac{\partial( \omega_N+ Z_1)}{\partial z}  \\
 \label{eq-48}
&&\ +f{k}\times (\omega_N+Z_1)+\nabla_H p_{b}-\frac{1}{\sqrt{C_0}}\int^{z}_{-1}\nabla_H (g_N+Z_2)dz'-\Delta\omega_N-\frac{\partial^2 \omega_N}{\partial z^2}=0,
\end{eqnarray}
\begin{eqnarray}\label{eq-70}
\frac{\partial g_N}{\partial t}+[(\omega_N+ Z_1)\cdot \nabla_H](g_N+ Z_2)+\Phi(\omega_N+Z_1)\frac{\partial (g_N+Z_2)}{\partial z}-\Delta g_N -\frac{\partial^2 g_N}{\partial z^2}=0.&
\end{eqnarray}
Taking the scalar product of (\ref{eq-48}) with $\omega_N$,  it follows that
\begin{eqnarray*}\notag
\frac{1}{2}\frac{d|\omega_N|^2}{dt}+\|\omega_N\|^2&=&-\Big(\omega_N, \big((\omega_N+Z_1)\cdot \nabla_H\big)(\omega_N+ Z_1)+\Phi(\omega_N+Z_1)\frac{\partial( \omega_N+ Z_1)}{\partial z}\Big)\nonumber\\
& &-\Big(\omega_N, f{k}\times (\omega_N+Z_1)\Big)-\Big(\omega_N, \nabla_H p_{b}-\frac{1}{\sqrt{C_0}}\int^{z}_{-1}\nabla_H (g_N+Z_2)dz'\Big).
\end{eqnarray*}
By integration by parts, we have
\[
\left(\omega_N, \Big((\omega_N+Z_1)\cdot \nabla_H \Big)\omega_N+\Phi(\omega_N+Z_1)\frac{\partial \omega_N}{\partial z}\right)=0, \quad (\omega_N, \nabla_H p_b)=0.
\]
Using Lemma \ref{le-2}, the H\"{o}lder inequality and the Young's inequality, we have
\begin{eqnarray*}
 \left|\left(\omega_N, \big((\omega_N+Z_1)\cdot \nabla_H\big)Z_1+\Phi(\omega_N+Z_1)\frac{\partial Z_1}{\partial z}\right)\right|&\leq& c |\omega_N|\|\omega_N+Z_1\|\|Z_1\|^{\frac{1}{2}}\|Z_1\|^{\frac{1}{2}}_2+c|\omega_N|\|\omega_N+Z_1\|\|Z_1\|^{\frac{1}{2}}_2\|Z_1\|^{\frac{1}{2}}_3\\
 &\leq&c\|Z_1\|_3\|\omega_N\|^2+\frac{1}{8}|\omega_N|^2+c\|Z_1\|^4_3+\frac{1}{8}\|\omega_N\|^2.
\end{eqnarray*}
 Since $(\omega_N, f{k}\times \omega_N)=0$, we get  $|(\omega_N, f{k}\times Z_1)|\leq\frac{1}{8}|\omega_N|^2+c|Z_1|^2$. Moreover, by the Young's inequality, we deduce that
 \begin{eqnarray*}
 \left|\left(\omega_N,\frac{1}{\sqrt{C_0}}\int^{z}_{-1}\nabla_H (g_N+Z_2)dz' \right)\right|&\leq& \frac{1}{2}\|\omega_N\|^2+\frac{1}{C_0}|g_N|^2+\frac{1}{C_0}|Z_2|^2.
\end{eqnarray*}
%
 Combining all inequalities in the above, we obtain
\begin{eqnarray}\notag
\frac{d|\omega_N|^2}{dt}+\frac{1}{2}\|\omega_N\|^2
&\leq& c\|Z_1\|_3\|\omega_N\|^2+c\|Z_1\|^4_3+\frac{1}{C_0}|g_N|^2+\frac{1}{C_0}|Z_2|^2+c|Z_1|^2\\
\label{eq-2-1}
&\leq& cM_8^{\frac{1}{2}}\|\omega_N\|^2+cM_8+\frac{2}{C_0}|g_N|^2.
\end{eqnarray}
Similarly, taking the scalar product on both sides of (\ref{eq-70}) with $g_N$, using H\"{o}lder inequality and Sobolev embedding theorem, it follows that
\begin{equation}\label{eq-2-2}
\frac{d|g_N|^2}{dt}+\frac{3}{2}\|g_N\|^2\leq cM_8|\omega_N|^2+cM^2_8|g_N|^2+2cM^2_8.
\end{equation}
Since $C_0\geq \frac{2}{\lambda_1}$, then $\frac{2}{C_0}|g_N|^2\leq \lambda_1|g_N|^2\leq \|g_N\|^2$.
When $M_8$ is sufficiently small, combining (\ref{eq-2-1}) and (\ref{eq-2-2}), we obtain
\begin{equation}\label{eq-49}
\frac{d|X_N|^2}{dt}+\frac{1}{4}\|X_N\|^2\leq cM_8 \quad {\rm on}\ (0,t).
\end{equation}
Applying Gronwall inequality, we get
\begin{equation}\notag
|X_N(t)|^2\leq {\rm{e}}^{-\frac{\mu_1}{4}t}|y_0|^2+\frac{cM_8}{\mu_1}.
\end{equation}
Then, we deduce from $|y_{0}|^2\leq 4\kappa$ that
\begin{equation}\notag
|X_N(t)|^2\leq 4\kappa {\rm{e}}^{-\frac{\mu_1}{4} t}+\frac{cM_8}{\mu_1}.
\end{equation}

Choosing $t$ sufficiently large and $M_8$ sufficiently small, we obtain (\ref{eq-47}), which yields Lemma \ref{le-6}. Indeed, when (\ref{eq-47}) holds,
\begin{eqnarray*}\notag
\mathbb{P}\left(|Y_N(t,y_{0})|^2\leq \delta_3\right)
&\geq & \mathbb{P}\left(|X_N(t,y_{0})|^2\leq \frac{\delta_3}{4},\quad  \sup_{(0,t)}\|Z(s)\|^2_3\leq M_8(\delta_3)\wedge \frac{\delta_3}{4}\right)\\ \notag
&\geq & \mathbb{P}\left(\sup_{(0,t)}\|Z(s)\|^2_3\leq M_8(\delta_3)\wedge \frac{\delta_3}{4}\right)\\ \notag
&\geq &  p_0\left(t,M_8(\delta_3)\wedge \frac{\delta_3}{4}\right )>0.
\end{eqnarray*}
Let $\tilde{M}_8(\delta_3)=M_8(\delta_3)\wedge \frac{\delta_3}{4}$
and $p_8(t, \delta_3)= p_0(t,\tilde{M}_8(\delta_3))$, then
 (\ref{eq-46}) holds. We complete the proof.
\end{proof}

As stated in Lemma \ref{le-6}, we have show the exponential moment estimates of the time entering into a small ball in $H$. In order to obtain such estimates for the space $\mathbb{H}_3$, three steps are needed.

\begin{lemma}\label{le-7}
Assume that  \textbf{Hypothesis H1} holds. Then, for any $\delta_4>0$, there exist $p_4(\delta_4)$ and $R_4(\delta_4)>0 $ such that for any $|y_0|^2\leq R_4$, we have for any $\Upsilon\leq 1$,
\[
\mathbb{P}\left(\|Y_N(\Upsilon, y_0)\|^2\leq \delta_4\right)\geq p_4.
\]
\end{lemma}
\begin{proof}
Using the decomposition $X_N=Y_N-P_N Z$ defined in Lemma \ref{le-5} and setting
\[
N(\omega)=\sup_{(0,\Upsilon)}\|Z(s,\omega)\|^2_3 \quad \rm{for} \ \omega\in \Omega.
\]
Let $\delta_4 >0$, assume that there exist $M_9(\delta_4)>0$, $R_4(\delta_4)>0$, such that for $\omega\in \Omega$,
\[
N(\omega)\leq M_9(\delta_4)\wedge \frac{\delta_4}{4} \  {\rm{implies}} \quad \|X_N(\Upsilon,\omega, y_0)\|^2\leq \frac{\delta_4}{4},
\]
provided $|y_0|^2\leq R_4(\delta_4)$. Then, we deduce Lemma \ref{le-7} holds from Lemma \ref{le-5} with  $M=M_9(\delta_4)\wedge\frac{\delta_4}{4}$.\\
Integrating (\ref{eq-49}), we obtain
\[
\frac{1}{4\Upsilon}\int^{\Upsilon}_{0}\|X_N(t)\|^2dt\leq \frac{1}{\Upsilon}|y_0|^2+cM_8,
\]
which yields,
\begin{equation}\label{eq-50}
\lambda\left(t\in(0,\Upsilon)\mid \|X_N(t)\|^2\leq \frac{8}{\Upsilon}|y_0|^2+8cM_8 \right)\geq \frac{\Upsilon}{2},
\end{equation}
where $\lambda$ denotes the Lebesgue measure on $(0, \Upsilon)$.
Set
\[
\tau_{\mathbb{H}_1}=\inf\Big\{t\in(0, \Upsilon)\mid \|X_N(t)\|^2\leq \frac{8}{\Upsilon}|y_0|^2+8cM_8\Big\}.
\]
We deduce from (\ref{eq-50}) and the continuity of $X_N$ that
\begin{equation}\label{eq-52}
\|X_N(\tau_{\mathbb{H}_1})\|^2\leq \frac{8}{\Upsilon}|y_0|^2+8cM_8.
\end{equation}
Taking inner product with $A_1\omega_N $ on both sides of (\ref{eq-48}) in $H$, we obtain
\begin{eqnarray}\notag
\frac{1}{2}\frac{d\|\omega_N\|^2}{dt}+\|\omega_N\|^2_2&=& -\left(A_1\omega_N,\left((\omega_N+Z_1)\cdot \nabla_H\right)(\omega_N+ Z_1)+\Phi(\omega_N+Z_1)\frac{\partial( \omega_N+ Z_1)}{\partial z}\right)\nonumber \\
& &-\left(A_1\omega_N, f{k}\times (\omega_N+Z_1)\right)-\left(A_1\omega_N, \nabla_H p_{b}-\frac{1}{\sqrt{C_0}}\int^{z}_{-1}\nabla_H (g_N+Z_2)dz'\right).\nonumber
\end{eqnarray}
Since
\[
|(A_1y, (x\cdot \nabla_H) z+(z\cdot \nabla_H) x)|\leq c\|y\|_2\|z\|^{\frac{1}{2}}\|z\|^{\frac{1}{2}}_2\|x\|,
\]
we have
\[
\begin{aligned}
&\left|\left(A_1\omega_N, \big((\omega_N+Z_1)\cdot \nabla_H\big)(\omega_N+ Z_1)+\Phi(\omega_N+Z_1)\frac{\partial( \omega_N+ Z_1)}{\partial z}\right)\right|\\
&\leq c\|\omega_N\|^{\frac{3}{2}}\|\omega_N\|^{\frac{3}{2}}_2+c\|Z_1\|_2\|\omega_N\|^2_2+c\|Z_1\|^2_2\|\omega_N\|_2+c\|\omega_N\|\|\omega_N\|^2_2+c\|Z_1\|_2\|\omega_N\|^2_2.
\end{aligned}
\]
By H\"{o}lder inequality and the Young's inequality, we have
\begin{eqnarray*}
\left|(A_1\omega_N, f{k}\times (\omega_N+Z_1)\right|&\leq& c|Z_1|\|\omega_N\|_2+C\|\omega_N\|^2,\\
\left|\left(A_1\omega_N, \nabla_H p_{b}-\frac{1}{\sqrt{C_0}}\int^{z}_{-1}\nabla_H (g_N+Z_2)dz'\right)\right|&\leq& \frac{1}{2}\|\omega_N\|^2_2+\frac{1}{C_0}\|g_N\|^2+\frac{1}{C_0}\|Z_2\|^2.
\end{eqnarray*}
 Thus, it follows that
\begin{equation}\label{eq-3-1}
\frac{d\|\omega_N\|^2}{dt}+\|\omega_N\|^2_2 \leq c\|\omega_N\|^6 +cM _9 +4cM^{\frac{1}{2}}_9\|\omega_N\|^2_2+\frac{2}{C_0}\|g_N\|^2+c\|\omega_N\|^2\|\omega_N\|^2_2.
\end{equation}
Similarly, taking the inner product with $A_2 g_N $ on both sides of (\ref{eq-70}) in $H_2$, we obtain
\begin{eqnarray}\label{eq-3-2}
\notag
\frac{d\|g_N\|^2}{dt}+\frac{3}{2}\|g_N\|^2_2 & \leq &  c\|g_N\|^2\|\omega_N\|^2\|\omega_N\|^2_2 +cM^2_9 +cM_9\|\omega_N\|^2_2\\
& &+cM^{\frac{1}{2}}_9\|g_N\|^2_2+c\|\omega_N\|^4\|g_N\|^2.
\end{eqnarray}
Since $C_0\geq \frac{2}{\lambda_1}$, then $\frac{2}{C_0}\|g_N\|^2\leq \lambda_1\|g_N\|^2\leq \|g_N\|^2_2$.
When $M_9$ is sufficiently small, combining (\ref{eq-3-1}) and (\ref{eq-3-2}),
\begin{equation}\notag
\frac{d\|X_N\|^2}{dt}+\|X_N\|^2_2\leq c\|\omega_N\|^6+c\|\omega_N\|^2\|\omega_N\|^2_2+c\|g_N\|^2\|\omega_N\|^2\|\omega_N\|^2_2+cM_9.
\end{equation}
Then
\[
\frac{d\|X_N\|^2}{dt}+\frac{1}{4}\|X_N\|^2_2+\left(\frac{1}{8}-c\|X_N\|^2-c\|X_N\|^4\right)\|\omega_N\|^2_2\leq c\|X_N\|^2(\|X_N\|^4-4K^2_1)+cM_9,
\]
 where  $K_1=\sqrt{\frac{\mu_1}{8c}}$ . Set
\[
\sigma_{\mathbb{H}_1}=\inf\left\{t\in(\tau_{\mathbb{H}_1},\Upsilon)\mid \|X_N(t)\|^2>2K_1\wedge \frac{1}{\sqrt{8c}}\right\}, \quad 2\tilde{K_1}=2K_1\wedge \frac{1}{\sqrt{8c}}.
\]
 Remark that on $(\tau_{\mathbb{H}_1},\sigma_{\mathbb{H}_1})$, we have
\begin{equation}\label{eq-55}
\frac{d\|X_N\|^2}{dt}+\frac{1}{4}\|X_N\|^2_2\leq cM_9.
\end{equation}
Integrating (\ref{eq-55}), we obtain
\begin{eqnarray}\notag
\|X_N(\sigma_{\mathbb{H}_1})\|^2+\frac{1}{4}\int^{\sigma_{\mathbb{H}_1}}_{\tau_{\mathbb{H}_1}}\|X_N(t)\|^2_2dt
&\leq& \|X_N(\tau_{\mathbb{H}_1})\|^2+cM_9(\sigma_{\mathbb{H}_1}-\tau_{\mathbb{H}_1})\\ \nonumber
&\leq& \|X_N(\tau_{\mathbb{H}_1})\|^2+cM_9\Upsilon\\
\label{eq-53}
&\leq& \|X_N(\tau_{\mathbb{H}_1})\|^2+cM_9.
\end{eqnarray}

From (\ref{eq-52}) and (\ref{eq-53}), we obtain that, for $M_8$, $M_9$ and $|y_0|^2$ sufficiently small,
\[
\|X_N(\sigma_{\mathbb{H}_1})\|^2 \leq \frac{\delta_4}{4}\wedge \tilde{K_1},
\]
which yields $\sigma_{\mathbb{H}_1}=\Upsilon$.  It follows that
\[
\|X_N(\Upsilon)\|^2 \leq \frac{\delta_4}{4},
\]
provided $M_8$, $M_9$ and $|y_0|^2$ sufficiently small.  Since
\begin{eqnarray*}\notag
\mathbb{P}\left(\|Y_N(\Upsilon,y_0)\|^2\leq \delta_4\right)
&\geq& \mathbb{P}\left(\|X_N(\Upsilon,y_0)\|^2\leq \frac{\delta_4}{4},\quad  \sup_{(0,\Upsilon)}\|Z(s)\|^2_3\leq M_9(\delta_4)\wedge\frac{\delta_4}{4}\right)\\ \notag
&\geq &\mathbb{P}\left(\sup_{(0,\Upsilon)}\|Z(s)\|^2_3\leq M_9(\delta_4)\wedge\frac{\delta_4}{4}\right)\\ \notag
&\geq & p_0\left(\Upsilon, M_9(\delta_4)\wedge\frac{\delta_4}{4}\right)>0.
\end{eqnarray*}
Let $\tilde{M}_9(\delta_4)=M_9(\delta_4)\wedge\frac{\delta_4}{4}$
and $p_4(\delta_4)=p_0(\Upsilon, \tilde{M}_9(\delta_4))$, we
obtain Lemma \ref{le-7}.
\end{proof}

\begin{lemma}\label{le-8}
Assume that  \textbf{Hypothesis H1} holds. Then, for any $\delta_5>0$, there exist $p_5(\delta_5)$ and $R_5(\delta_5)>0 $ such that for any $y_0$ verifying $\|y_0\|^2\leq R_5$, we have for any $\Upsilon\leq 1$,
\[
\mathbb{P}\left(\|Y_N(\Upsilon,y_0)\|^2_2\leq \delta_5\right)\geq p_5.
\]
\end{lemma}

 \begin{proof}
When $M_9$ sufficiently small  and $\|y_0\|^2+cM_9\leq \tilde{K}_1$, we have $\tau_{\mathbb{H}_1}=0$ and $\sigma_{\mathbb{H}_1}=\Upsilon$, then it follows from (\ref{eq-55}) that
\[
\frac{1}{4}\int^{\Upsilon}_{0}\|X_N(t)\|^2_2dt\leq \|y_0\|^2+cM_9.
\]
Similar to the proof of (\ref{eq-52}), we deduce that there exists a stopping time $\tau_{\mathbb{H}_2}\in (0,\Upsilon)$ such that
\begin{equation}\label{eq-56}
\|X_N(\tau_{\mathbb{H}_2})\|^2_2\leq\frac{8}{\Upsilon}(\|y_0\|^2+cM_9),
\end{equation}
provided $M_9$ and $\|y_0\|^2$ are sufficiently small.

Using the same method as Lemma \ref{le-7}, let $\delta_5 >0$, assume that there exist $M_{10}(\delta_5)>0$ and $R_5(\delta_5)>0$  such that for $\omega\in \Omega$,
\[
N(\omega)\leq M_{10} (\delta_5)\wedge \frac{\delta_5}{4} \  {\rm{implies}} \quad \|X_N(\Upsilon,\omega, y_0)\|^2_2\leq \frac{\delta_5}{4},
\]
provided $\|y_0\|^2\leq R_5(\delta_5)$. Then, we conclude Lemma \ref{le-8} holds from Lemma \ref{le-5} with
$M=M_{10}(\delta_5)\wedge\frac{\delta_5}{4}$.

Taking the scalar product with $A^2_1 \omega_N$ on both sides of (\ref{eq-48}), we obtain
\begin{eqnarray}\notag
\frac{1}{2}\frac{d\|\omega_N\|^2_2}{dt}+\|\omega_N\|^2_3&=&-\left(A^2_1\omega_N, \Big((\omega_N+Z_1)\cdot \nabla_H\Big)(\omega_N+ Z_1)+\Phi(\omega_N+Z_1)\frac{\partial( \omega_N+ Z_1)}{\partial z}\right)\notag \\
& &-\left(A^2_1\omega_N, f{k}\times (\omega_N+Z_1)\right)-\left(A^2_1\omega_N, \nabla_H p_{b}-\frac{1}{\sqrt{C_0}}\int^{z}_{-1}\nabla_H (g_N+Z_2)dz'\right).\notag
\end{eqnarray}
Since
\[
|(A^2_1y, (x\cdot \nabla_H) z+(z\cdot \nabla_H) x)|\leq c\|y\|_3\|z\|_2\|x\|_2\leq \varepsilon\|y\|^2_3+c(\|z\|^4_2+\|x\|^4_2),
\]
we obtain
\[
\begin{aligned}
&\left|\left(A^2_1\omega_N,\Big((\omega_N+Z_1)\cdot \nabla_H\Big)(\omega_N+ Z_1)+\Phi(\omega_N+Z_1)\frac{\partial( \omega_N+ Z_1)}{\partial z}\right)\right|\\
&\leq c\|\omega_N\|^{\frac{3}{2}}_3\|\omega_N\|^{\frac{1}{2}}_2\|\omega_N\|+c\|\omega_N\|_3\|\omega_N\|^{\frac{3}{2}}_2\|\omega_N\|^{\frac{1}{2}}
+c\|\omega_N\|_3\|\omega_N\|^{\frac{1}{2}}_2\|\omega_N\|^{\frac{1}{2}}\|Z_1\|_2\\
& \quad +c\|\omega_N\|^2_3\|Z_1\|_3
+c\|\omega_N\|_3\|Z_1\|^2_3+c\|\omega_N\|^{\frac{3}{2}}_3\|\omega_N\|^{\frac{3}{2}}_2
+c\|Z_1\|_2\|\omega_N\|^2_2\\
& \quad +c\|Z_1\|^2_2\|\omega_N\|_2+c\|\omega_N\|\|\omega_N\|^2_2+c\|Z_1\|_2\|\omega_N\|^2_2\\
&\leq  \varepsilon\|\omega_N\|^2_3+c\|\omega_N\|^6_2+c\|Z_1\|^2_3\|\omega_N\|^2_2+c\|Z_1\|^4_3.
\end{aligned}
\]
By H\"{o}lder inequality and the Young' inequality, we have
\[
|(A^2_1\omega_N, Z_1)|\leq c\|Z_1\|\|\omega_N\|_3,
\]
\[
\left|\left(A^2_1\omega_N, \nabla_H p_{b}-\frac{1}{\sqrt{C_0}}\int^{z}_{-1}\nabla_H (g_N+Z_2)dz'\right)\right|\leq \frac{1}{2}\|\omega_N\|^2_3+\frac{1}{C_0}\|g_N\|^2_2+\frac{1}{C_0}\|Z_2\|^2_2.
\]
Combining the above inequalities, applying H\"{o}lder inequality and the Young' inequality,  we obtain
\begin{equation}\label{eq-4-1}
\frac{d\|\omega_N\|^2_2}{dt}+\|\omega_N\|^2_3\leq c\|\omega_N\|^6_2+c\|\omega_N\|^4_2+cM^{\frac{1}{2}}_{10}\|\omega_N\|^2_3+\frac{2}{C_0}\|g_N\|^2_2+cM^2_{10}+\varepsilon \|\omega_N\|^2_3.
\end{equation}
Similarly, taking the scalar product with $A^2_2 g_N$ on both sides of (\ref{eq-70}), we obtain
\begin{eqnarray}\notag
\frac{d\|g_N\|^2_2}{dt}+2\|g_N\|^2_3 & \leq & \varepsilon \|g_N\|^2_3+c\|g_N\|^2_2\|\omega_N\|^4_2+c\|g_N\|^2_2\|\omega_N\|^2_2\\
\label{eq-4-2}
&  &+cM_{10}\|\omega_N\|^2_3+cM_{10}\|g_N\|^2_3+cM^2_{10}.
\end{eqnarray}

Choosing $M_{10}$ small enough, since $\frac{2}{C_0}\|g_N\|^2_2\leq \lambda_1\|g_N\|^2_2\leq\|g_N\|^2_3$,  by (\ref{eq-4-1}) and (\ref{eq-4-2}), we obtain
\begin{equation}\label{eq-58}
\frac{d\|X_N\|^2_2}{dt}+\frac{1}{4}\|X_N\|^2_3\leq c\|X_N\|^2_2\left(\|X_N\|^2_2+\|X_N\|^4_2-4K^2_2\right)+cM_{10}.
\end{equation}
where $K_2$ is defined similar to the above.
Setting $\sigma_{\mathbb{H}_2}=\inf\{t\in (\tau_{\mathbb{H}_2}, \Upsilon)\mid \|X_N\|^2_2\geq 2K_2\}$, integrating (\ref{eq-58}) on $(\tau_{\mathbb{H}_2},\sigma_{\mathbb{H}_2})$, we obtain
\begin{equation}\label{eq-57}
\|X_N(\sigma_{\mathbb{H}_2})\|^2_2+\frac{1}{4}\int^{\sigma_{\mathbb{H}_2}}_{\tau_{\mathbb{H}_2}}\|X_N(t)\|^2_3dt\leq \|X_N(\tau_{\mathbb{H}_2})\|^2_2+cM_{10}.
\end{equation}

Combining (\ref{eq-56}) and (\ref{eq-57}). We obtain that, for $M_{9}$, $M_{10}$ and $\|y_0\|^2$\ sufficiently small,
\[
\|X_N(\sigma_{\mathbb{H}_2})\|^2_2\leq \frac{\delta_5}{4}\wedge K_2.
\]
It follows that $\sigma_{\mathbb{H}_2}=\Upsilon$ and
\[
\|X_N(\Upsilon)\|^2_2\leq \frac{\delta_5}{4},
\]
provided $M_{9}$, $M_{10}$ and $\|y_0\|^2$ sufficiently small. Since
\begin{eqnarray}\notag
\mathbb{P}\left(\|Y_N(\Upsilon,y_0)\|^2_2\leq \delta_5\right)
&\geq & \mathbb{P}\left(\|X_N(\Upsilon,y_0)\|^2\leq \frac{\delta_5}{4},\quad \sup_{(0,\Upsilon)}\|Z(s)\|^2_3\leq M_{10}(\delta_5)\wedge\frac{\delta_5}{4}\right)\\ \notag
&\geq & \mathbb{P}\left(\sup_{(0,\Upsilon)}\|Z(s)\|^2_3\leq M_{10}(\delta_5)\wedge\frac{\delta_5}{4}\right)\\ \notag
&\geq  & p_0\left(\Upsilon, M_{10}(\delta_5)\wedge\frac{\delta_5}{4}\right)>0, \notag
\end{eqnarray}
let
$\tilde{M}_{10}(\delta_5)=M_{10}(\delta_5)\wedge\frac{\delta_5}{4}$
and  $p_5(\delta_5)=p_0(\Upsilon, \tilde{M}_{10}(\delta_5))$,
we obtain Lemma \ref{le-8}.
\end{proof}

\begin{lemma}\label{le-9}
Assume that  \textbf{Hypothesis H1} holds. Then, for any $\delta_6>0$, there exist $p_6(\delta_6)$ and $R_6(\delta_6)>0 $ such that for any $y_0$ verifying $\|y_0\|^2_2\leq R_6$, we have for any $\Upsilon\leq 1$,
\[
\mathbb{P}\left(\|Y_N(\Upsilon,y_0)\|^2_3\leq \delta_6\right)\geq p_6.
\]
\end{lemma}

 \begin{proof}
 When $M_{10}$ sufficiently small and $\|y_0\|^2_2+cM_{10}\leq K_2$, we have $\tau_{\mathbb{H}_2}=0$ and $\sigma_{\mathbb{H}_2}=\Upsilon$. Taking into account (\ref{eq-58}), we obtain
\[
\frac{1}{4}\int^{\Upsilon}_{0}\|X_N(t)\|^2_3dt\leq \|y_0\|^2_2+cM_{10}.
\]
Similar to the proof of (\ref{eq-52}), we know that there exists a stopping time $\tau_{\mathbb{H}_3}\in (0,\Upsilon)$ such that
\begin{equation}\label{eq8}
\|X_N(\tau_{\mathbb{H}_3})\|^2_3\leq \frac{8}{\Upsilon}(\|y_0\|^2_2+cM_{10}).
\end{equation}
Similar to Lemma \ref{le-7}, let $\delta_6 >0$, assume that there exist $M_{11}(\delta_6)>0$ and $R_6(\delta_6)>0$ such that for $\omega\in \Omega$,
\[
N(\omega)\leq M_{11}(\delta_6)\wedge \frac{\delta_6}{4} \  {\rm{implies}} \quad \|X_N(\Upsilon, \omega, y_0)\|^2_3\leq \frac{\delta_6}{4},
\]
provided $\|y_0\|^2_2\leq R_6(\delta_6)$. Then, we deduce that Lemma \ref{le-9} holds from Lemma \ref{le-5} with
$M=M_{11}(\delta_6)\wedge\frac{\delta_6}{4}$.

Taking the scalar product with $A^3_1 \omega_N$ on both sides of (\ref{eq-48}), we obtain
\begin{eqnarray}\notag
\frac{1}{2}\frac{d\|\omega_N\|^2_3}{dt}+\|\omega_N\|^2_4&=& -\left(A^3_1\omega_N, \Big((\omega_N+Z_1)\cdot \nabla_H\Big)(\omega_N+ Z_1)+\Phi(\omega_N+Z_1)\frac{\partial( \omega_N+ Z_1)}{\partial z}\right) \notag \\
& &-\left(A^3_1\omega_N, f{k}\times (\omega_N+Z_1)\right)-\left(A^3_1\omega_N, \nabla_H p_{b}-\frac{1}{\sqrt{C_0}}\int^{z}_{-1}\nabla_H (g_N+Z_2)dz'\right),\notag
\end{eqnarray}
Since
\[
|(A^3_1y, (x\cdot \nabla_H) z+(z\cdot \nabla_H) x)|\leq c\|y\|_4\|z\|_3\|x\|_3\leq \varepsilon\|y\|^2_4+c(\|z\|^4_3+\|x\|^4_3),
\]
we obtain
\[
\begin{aligned}
&\left|\left(A^3_1\omega_N, \big((\omega_N+Z_1)\cdot \nabla_H\big)(\omega_N+ Z_1)+\Phi(\omega_N+Z_1)\frac{\partial( \omega_N+ Z_1)}{\partial z}\right)\right|\\
&\leq c\|\omega_N\|_4\|\omega_N\|_3\|\omega_N\|_2+c\|\omega_N\|^2_4\|Z_1\|_3+c\|\omega_N\|^2_3\|Z_1\|^2_3\\
&\quad +c\|\omega_N\|_4\|Z_1\|^2_3+c\|\omega_N\|_4\|\omega_N\|^2_3+c\|\omega_N\|^{\frac{3}{2}}_4\|\omega_N\|^{\frac{1}{2}}_3\|Z_1\|_2\\
&\leq \varepsilon\|\omega_N\|^2_4+c\|\omega_N\|^4_3+c\|Z_1\|^2_3\|\omega_N\|^2_4+c\|Z_1\|^4_3+c\|Z_1\|^4_3\|\omega_N\|^2_3.
\end{aligned}
\]
By H\"{o}lder inequality and the Young' inequality, we have
\[
|(A^3_1\omega_N, Z_1)|\leq c\|Z_1\|_2\|\omega_N\|_4,
\]
\[
\left|\left(A^3_1\omega_N, \nabla_H p_b-\frac{1}{\sqrt{C_0}}\int^{z}_{-1}\nabla_H (g_N+Z_2)dz'\right)\right|\leq \frac{1}{2}\|\omega_N\|^2_4+\frac{1}{C_0}\|g_N\|^2_3+\frac{1}{C_0}\|Z_2\|^2_3.
\]
Combining the above inequalities, we obtain
\begin{equation}\label{eq-5-1}
\frac{d\|\omega_N\|^2_3}{dt}+\|\omega_N\|^2_4\leq \varepsilon \|\omega_N\|^2_4+ c\|\omega_N\|^4_3+cM_{11}+cM_{11}\|\omega_N\|^2_3+\frac{2}{C_0}\|g_N\|^2_3.
\end{equation}
Similarly, taking the scalar product with $A^3_2 g_N$ on both sides of (\ref{eq-70}), we obtain
\begin{equation}\label{eq-5-2}
\frac{d\|g_N\|^2_3}{dt}+2\|g_N\|^2_4\leq \varepsilon \|g_N\|^2_4+ c\|\omega_N\|^2_3\|g_N\|^2_3+cM_{11}\|g_N\|^2_4+cM_{11}\|\omega_N\|^2_4+cM^2_{11}.
\end{equation}

Choosing $M_{11}$ sufficiently small, since $\frac{2}{C_0}\|g_N\|^2_3\leq \lambda_1\|g_N\|^2_3\leq \|g_N\|^2_4$, by (\ref{eq-5-1}) and (\ref{eq-5-2}), we obtain
\[
\frac{d\|X_N\|^2_3}{dt}+\frac{3}{4}\|X_N\|^2_4\leq c\|X_N\|^4_3+cM_{11},
\]
then
\begin{equation}\label{eq6}
\frac{d\|X_N\|^2_3}{dt}+\frac{1}{4}\|X_N\|^2_4\leq c\|X_N\|^2_3\left(\|X_N\|^2_3-2K_3\right)+cM_{11},
\end{equation}
where $K_3$ is defined similarly to the above.
Setting $\sigma_{\mathbb{H}_3}=\inf\{t\in (\tau_{\mathbb{H}_3}, \Upsilon)\mid \|X_N(t)\|^2_3\geq {2K_3}\}$ . Integrating (\ref{eq6}) on $(\tau_{\mathbb{H}_3}, \sigma_{\mathbb{H}_3})$, we obtain
\[
\|X_N(\sigma_{\mathbb{H}_3})\|^2_3+\frac{1}{4}\int^{\sigma_{\mathbb{H}_3}}_{\tau_{\mathbb{H}_3}}\|X_N(t)\|^2_4dt\leq \|X_N(\tau_{\mathbb{H}_3})\|^2_3+cM_{11}.
\]
Taking (\ref{eq8}) into account and choosing $M_{10}$, $M_{11}$ and $\|y_0\|^2_2$  sufficiently small, we obtain
\[
\|X_N(\sigma_{\mathbb{H}_3})\|^2_3\leq \frac{\delta_6}{4}\wedge K_3.
\]
It follows that $\sigma_{\mathbb{H}_3}=\Upsilon$ and that
\[
\|X_N(\Upsilon)\|^2_3\leq \frac{\delta_6}{4},
\]
provided $M_{10}$, $M_{11}$ and $\|y_0\|^2_2$ sufficiently small. Since
\begin{eqnarray}\notag
\mathbb{P}\left(\|Y_N(\Upsilon,y_0)\|^2_2\leq \delta_6\right)
&\geq& \mathbb{P}\left(\|X_N(\Upsilon,y_0)\|^2\leq \frac{\delta_6}{4},\quad \sup_{(0,\Upsilon)}\|Z(s)\|^2_3\leq M_{11}(\delta_6)\wedge\frac{\delta_6}{4}\right)\\ \notag
&\geq& \mathbb{P}\left(\sup_{(0,\Upsilon)}\|Z(s)\|^2_3\leq M_{11}(\delta_6)\wedge\frac{\delta_6}{4}\right)\\ \notag
&\geq & p_0\left(\Upsilon, M_{11}(\delta_6)\wedge\frac{\delta_6}{4}\right)>0, \notag
\end{eqnarray}
let
$\tilde{M}_{11}(\delta_6)=M_{11}(\delta_6)\wedge\frac{\delta_6}{4}$
and $p_6(\delta_6)=p_0(\Upsilon, \tilde{M}_{11}(\delta_6))$,
we obtain Lemma \ref{le-9}.
\end{proof}

Now, we are ready to prove Proposition \ref{pr-2}.

\begin{flushleft}
\textbf{Proof of
Proposition \ref{pr-2}.}\quad
Set
\[
\begin{aligned}
&\delta_6= \delta, \ \delta_5= R_6(\delta_6), \delta_4= R_5(\delta_5), \delta_3= R_4(\delta_4),
\\
&p_4=p_4(\delta_4),p_5=p_5(\delta_5),p_6=p_6(\delta_6),p_1=(p_4p_5p_6)^2.
\end{aligned}
\]
By the definition of $\tau_{L^2}$, we have
\[
|Y^1(\tau_{L^2})|^2\vee|Y^2(\tau_{L^2})|^2\leq R_4(\delta_4).
\]
In the following, we will prove this proposition in three cases.

\textbf{The first case:} $\|Y^1(\tau_{L^2})\|^2_3\vee\|Y^2(\tau_{L^2})\|^2_3\leq \delta$, which obviously yields
\begin{equation}\label{eq-41}
\mathbb{P}\left(\min_{k=0,1,2,3}\max_{i=1,2}\|Y^{i}(\tau_{L^2}+k\Upsilon)\|^2_3\leq\delta \mid \left(Y^1(\tau_{L^2}),Y^2(\tau_{L^2})\right)\right)\geq p_1.
\end{equation}

\textbf{The second case:} $Y^1(\tau_{L^2})=Y^2(\tau_{L^2})=y_0$ with $\|y_0\|^2_3>\delta$.
Combining Lemmas \ref{le-7}-\ref{le-9}, we deduce from the strong Markov property of $Y_N$ that
\[
\mathbb{P}(\|Y_N(3\Upsilon, y_0)\|^2_3\leq \delta)\geq p_6p_5p_4,
\]
provided $|y_0|^2\leq R_4$. In that case, $Y^1(\tau_{L^2}+3\Upsilon)=Y^2(\tau_{L^2}+3\Upsilon)$. Hence, since the law of $Y^1(\tau_{L^2}+3\Upsilon)$ conditioned by $(Y^1(\tau_{L^2}),Y^2(\tau_{L^2}))= (y^1_0, y^2_0)$  is $\mathcal{D}(Y_N(3\Upsilon, y_0))$, it follows that
\[
\mathbb{P}\left(\max_{i=1,2}\|Y^i(\tau_{L^2}+3\Upsilon)\|^2_3\leq \delta \mid \left(Y^1(\tau_{L^2}),Y^2(\tau_{L^2})\right)\right)\geq p_6p_5p_4\geq p_1,
\]
which implies (\ref{eq-41}) holds.

\textbf{The third case:} $Y^1(\tau_{L^2})\neq Y^2(\tau_{L^2})$ and $\|Y^1(\tau_{L^2})\|^2_3 \vee \|Y^2(\tau_{L^2})\|^2_3>\delta$. In that case, $\big(Y^1(\tau_{L^2}+\Upsilon),Y^2(\tau_{L^2}+\Upsilon)\big)$ conditioned by $\big(Y^1(\tau_{L^2}),Y^2(\tau_{L^2})\big)$ are independent. Since the law of $Y^i(\tau_{L^2}+\Upsilon)$ conditioned by $\big(Y^1(\tau_{L^2}),Y^2(\tau_{L^2})\big)=(y^1_0,y^2_0)$ is $\mathcal{D}(Y_N(\Upsilon, y^i_0))$, it follows from Lemma \ref{le-7} that
\begin{eqnarray}\label{eq-65}
\mathbb{P}\left(\max_{i=1,2}\|Y^i(\tau_{L^2}+\Upsilon)\|^2\leq \delta_4\mid \left(Y^1(\tau_{L^2}),Y^2(\tau_{L^2})\right)\right)\geq p^2_4.
\end{eqnarray}
Then, we distinguish three cases for $(Y^i(\tau_{L^2}+\Upsilon))_{i=1,2}$ similar to $(Y^i(\tau_{L^2}))_{i=1,2}$: $\|Y^1(\tau_{L^2}+\Upsilon)\|^2_3\vee\|Y^2(\tau_{L^2}+\Upsilon)\|^2_3\leq \delta$, $Y^1(\tau_{L^2}+\Upsilon)=Y^2(\tau_{L^2}+\Upsilon)=y_1$ with $\|y_1\|^2_3>\delta$, $Y^1(\tau_{L^2}+\Upsilon)\neq Y^2(\tau_{L^2}+\Upsilon)$ and $\|Y^1(\tau_{L^2}+\Upsilon)\|^2_3 \vee \|Y^2(\tau_{L^2}+\Upsilon)\|^2_3>\delta$. For the front two cases, using the method similar to the first and second case, respectively, and combining (\ref{eq-65}), we have (\ref{eq-41}) holds. For the last case, we know that $\big(Y^1(\tau_{L^2}+2\Upsilon),Y^2(\tau_{L^2}+2\Upsilon)\big)$ conditioned by $\big(Y^1(\tau_{L^2}+\Upsilon),Y^2(\tau_{L^2}+\Upsilon)\big)$ are independent. By Lemma \ref{le-8},  we have
\begin{eqnarray}\label{eq-66}
\mathbb{P}\left(\min_{k=1,2}\max_{i=1,2}\|Y^i(\tau_{L^2}+k\Upsilon)\|^2_2\leq \delta_5\mid \left(Y^1(\tau_{L^2}+\Upsilon),Y^2(\tau_{L^2}+\Upsilon)\right)\right)\geq p^2_5,
\end{eqnarray}
provided
\[
\max_{i=1,2}\|Y^i(\tau_{L^2}+\Upsilon)\|^2\leq \delta_4.
\]

Now, we distinguish three cases for $(Y^i(\tau_{L^2}+2\Upsilon))_{i=1,2}$ similarly to the above: $\|Y^1(\tau_{L^2}+2\Upsilon)\|^2_3\vee\|Y^2(\tau_{L^2}+2\Upsilon)\|^2_3\leq \delta$, $Y^1(\tau_{L^2}+2\Upsilon)=Y^2(\tau_{L^2}+2\Upsilon)=y_2$ with $\|y_2\|^2_3>\delta$, $Y^1(\tau_{L^2}+2\Upsilon)\neq Y^2(\tau_{L^2}+2\Upsilon)$ and $\|Y^1(\tau_{L^2}+2\Upsilon)\|^2_3 \vee \|Y^2(\tau_{L^2}+2\Upsilon)\|^2_3>\delta$. For the front two cases, we also can obtain (\ref{eq-41}). For the last case, $\big(Y^1(\tau_{L^2}+3\Upsilon),Y^2(\tau_{L^2}+3\Upsilon)\big)$ conditioned by $\big(Y^1(\tau_{L^2}+2\Upsilon),Y^2(\tau_{L^2}+2\Upsilon)\big)$ are independent. By Lemma  \ref{le-9}, we have
\begin{eqnarray}\label{eq-67}
\mathbb{P}\left(\min_{k=2,3}\max_{i=1,2}\|Y^{i}(\tau_{L^2}+k\Upsilon)\|^2_3\leq\delta \mid \left(Y^1(\tau_{L^2}+2\Upsilon),Y^2(\tau_{L^2}+2\Upsilon)\right)\right)\geq p^2_6,
\end{eqnarray}
provided
\[
\max_{i=1,2}\|Y^i(\tau_{L^2}+2\Upsilon)\|^2_2\leq \delta_5.
\]
Combining (\ref{eq-65})-(\ref{eq-67}), we deduce (\ref{eq-41}) holds for the last case.

Thus, we have proved (\ref{eq-41}) is  true almost surely. Integrating (\ref{eq-41}), we obtain
\begin{equation}\label{equa-2}
\mathbb{P}\left(\min_{k=0,1,2,3}\max_{i=1,2}\|Y^{i}(\tau_{L^2}+k\Upsilon)\|^2_3\leq\delta \right)\geq p_1.
\end{equation}
Combining Lemma  \ref{le-6} and (\ref{equa-2}), we conclude the result.
\end{flushleft}

$\hfill\blacksquare$




\noindent{\small {\bf  Acknowledgements}\ \  This work is partly supported by National Natural Science Foundation of China (No. 11431014, 11401557, 11801032). Key Laboratory of Random Complex Structures and Data Science, Academy of Mathematics and Systems Science, Chinese Academy of Sciences (No. 2008DP173182). China Postdoctoral Science Foundation funded project (No. 2018M641204).


\def\refname{ References}

\end{document}